\documentclass[aos]{imsart}

\startlocaldefs

\usepackage{graphicx, amssymb, amsmath, hyperref, bm, color, tikz, fullpage}
\usepackage[boxed]{algorithm2e}
\renewcommand{\bar}{\overline}

\renewcommand{\hat}{\widehat}

\begin{document}

\title{Testing for Tail Behavior Using Extreme Spacings}%
\author{Javier Rojo*}%
\author{ and Richard C. Ott\\}%

{\it Rice University and Mesa State College}
\address{
$
\begin{array}{llll}
\hspace*{-2.0in} $Department of Statistics$ & \;\;\;\;\;& & $Department of Computer Science, Mathematics$\\
\hspace*{-2.0in} $Rice University$ &\hspace*{1.0in}\;\;\;& \;\;\;\;\; &	$and Statistics$ \\
\hspace*{-2.0in} $MS 138$   &	\;\;\;\;\;& &	$Mesa State College$ \\
\hspace*{-2.0in}$6100 Main St. $ &\;\;\;\;\  & & $1100 North Avenue$ \\
 \hspace*{-2.0in}$Houston, Texas 77005$ & & & $Grand Junction, CO 81501$\\
 \hspace*{-2.0in}$Email: jrojo@rice.edu$ & & & $Email: rott@mesastate.edu$\\
 %$          $ & & & \\
 \end{array}$}

\let\oldthefootnote\thefootnote
\renewcommand{\thefootnote}{\fnsymbol{footnote}}
\footnotetext[1]{Research supported by The National Science Foundation grants
DMS-0532346 and DMS-0552590, and National Security Agency award H98230-06-1-0099. 
}
\let\thefootnote\oldthefootnote

\begin{abstract}
\small{Methodologies to test
hypotheses about the tail-heaviness of an underlying distribution 
are introduced based on %the 
results of Rojo (1996)\nocite{Rojo96} 
using the limiting behavior of the extreme spacings. The tests are consistent and have point-wise
robust levels in the sense of Lehmann (2005)\nocite{Leh05} and Lehmann and Loh (1990)\nocite{Leh90}. Simulation results based on these new methodologies indicate that the tests 
exhibit good control of the probability of Type I error and have good power 
properties for finite sample sizes. The tests are compared with a test proposed by Bryson (1974)\nocite{Bryson74} and it is seen that, although Bryson's test is competitive with the tests proposed here, Bryson's test does not have point-wise robust levels.
The operating characteristics of the tests are also explored when
the data is blocked.  It turns out that the power increases substantially by 
blocking. The methodology is illustrated by analyzing various data sets.}
\vspace*{.1in}

\end{abstract}
\pagenumbering{arabic}
\maketitle
%---------------------------------
\section    {\large Introduction}
%----------------------------------
Tail behavior of a probability distribution plays an important role in
various applications including hydrology,
aerospace engineering, meteorology, insurance, and finance. Lehmann (1988)\nocite{Lehmann} proposed
a pure-tail ordering in connection with the comparison, in terms of efficiency,
of location experiments.
In classical
extreme value theory, the Extremal Types Theorem (ETT), (the Three Type 
Theorem as it is also known), classifies
the right tail of a distribution according to the asymptotic distribution of the
standardized maximum. Thus, the distribution $F$ is short-, medium-, or long-tailed
depending on whether $F$ is in the domain of attraction of the Weibull, Gumbel,
or Fr\'{e}chet distributions. It is well-known, however, that the
 limiting distribution for the standardized maximum does not exist for 
 all distributions.
For instance, any distribution that assigns positive mass to the right
endpoint of its support cannot be classified by the ETT. Another possible
drawback of the classical categorization of probability laws using the
ETT is that
the class of medium-tailed distributions may be too large. For
instance, Schuster (1984)\nocite{Sch84} argues that, ''the statistician
considers the normal distribution shorter than the exponential
 which is in turn shorter than the lognormal
distribution''. Yet all three are in the domain of attraction of
the Gumbel distribution. Thus, there is a need to classify
distributions by alternative schemes.  

Another possible avenue for such classification may be
obtained through the tail-heaviness of a distribution. 
There exists a robust literature on orderings that attempt to
order distributions according to tail-heaviness. Some of the early
work allows for the middle part of the
distribution to affect the tail ordering. See, e.g.,  Loh (1984)\nocite{Loh}, 
Doksum (1969)\nocite{Doksum}, and Lehmann (1988)\nocite{Lehmann}.
By contrast, in a series of papers, Rojo
(1988, 1992, 1993, 1996)\nocite{Rojo88}\nocite{Rojo92}\nocite{Rojo93}\nocite{Rojo96}
proposes pure-tail orderings which allow for pure-tail comparisons without many 
of the technical assumptions required by other approaches. 

Various methods for distinguishing between exponential and power-tails have been proposed 
in the literature and are used in practice. The more popular ones are based on plotting various
quantities whose behavior depend on the tail-behavior of the underlying distribution.
For example, Heyde and Kou (2004)\nocite{Heyde04} discuss methods based on plotting 
the mean residual life function, Q-Q plots, conditional moment generating functions, Hill estimation, 
and likelihood methods. Heyde and Kou (2004)\nocite{Heyde04} argue that these methods are qualitative without any support for their "statistical precision". The Hill estimator is popular in practice
but it has its share of problems, including its undesirable behavior as represented by the "Hill's horror
plots" as illustrated in Embrechts {\it et al} (1997)\nocite{Embrechts}, in spite of the various results 
concerning its asymptotic properties; moreover, its use should be restricted to the case of power-tails 
since the Hill estimation may be misleading, as it does not provide alerts, when applied to other types of tails. Thus the use of Hill estimation may be inappropriate for distinguishing power tails from 
other classes of tails. Somewhat surprisingly, Heyde and Kou (2004) conclude that it may be necessary to have sample sizes in the tens,  and sometimes in the hundreds, of thousands to be able to differentiate between power and exponential tails. The main reason for this is that the large quantiles of exponentially-tailed distributions may actually exceed the counterpart quantiles of power-like tails. This characteristic will be observed in our simulation work. One way to ameliorate this problem is by blocking the data.

The purpose of this paper is to develop methodologies to test
hypotheses about the tail-heaviness of a distribution based on the
results of Rojo (1996)\nocite{Rojo96}. The paper will focus on the right tail of 
the distribution, but analogous results are easily seen to hold for the left
tail by considering instead the behavior of random variables $\{-X_i, \dots, i=1,\dots,n\}$. 
When the underlying distribution is symmetric about zero, one may take $\{|X_i|, \dots, i=1,\dots,n\}$ in effect doubling the sample size. Theorem 3.1 and Corollary 4.2 in
Rojo (1996)\nocite{Rojo96}  provide the results needed to develop
methodologies to test the hypothesis that data arises from a
medium-tailed distribution against an alternative of a short- or
long-tailed distribution based on the asymptotic distribution of the 
extreme spacing $X_{(n)}-X_{(n-1)}$, where $X_{(k)}$ represents the $k^{th}$ order statistic
from a random sample $X_{1}, \dots,X_{n}$ from $F$. The distribution $F$ is assumed to be continuous and strictly increasing throughout this work.
 
The organization of the paper is as follows: Sections 1 and 2 provide the introductory 
material and a brief discussion of classification schemes developed by Parzen (1979)\nocite{Parzen}
and Schuster (1984)\nocite{Sch84}. Section 3 discusses the most relevant results
from Rojo (1996)\nocite{Rojo96} and section 4 discusses a new test for tail-heaviness.
The test is consistent against short- and long-tailed alternatives and the level of the test is point-wise robust. Simulation results indicate that for small sample sizes the test exhibits good control of
the probability of Type I error, and has good power properties. A comparison with a test proposed 
by Bryson (1974) concludes that, although Bryson's test behaves well against distributions with linear mean residual life functions, its power is not good against distributions with quadratic mean residual 
life functions and its probability of type I error is close to 1 for the gamma and log-gamma distributions (which are medium-tailed) and hence may not be a good choice for the testing situation of interest in this work. A simulation study
is discussed where that data is blocked  to increase the power of the test. Finally, the methodology is illustrated by applying it to several published data sets: maximum discharge of the Feather river, glass breaking strength, and the Belgian Secure Re claim size data.

%%%%%%%%%%%%%%%%%%%%%%%%%%%%%%%%%%%%%%%%%%%%%%
\section {\large Classification Based on the Density-quantile Function}

%%%%%%%%%%%%%%%%%%%%%%%%%%%%%%%%%%%%%%%%%%%%%%
Parzen (1979)\nocite{Parzen} argued that many 
distributions have density-quantile functions of the form, 
\begin{equation}\label{denQ}
fQ(u)\backsim
(1-u)^{\alpha},\  \alpha >0,\end{equation} 
where $f$ denotes the density function, $Q$ is the quantile function
(the left-continuous inverse of $F$), and $g_{1}(u)\backsim%
g_{2}(u) $ means that $g_{1}(u)/g_{2}(u)$ tends to a positive finite constant
as $u\rightarrow 1.$ The parameter $\alpha$ is called the tail
exponent, and Parzen (1979) defined a distribution to
be short-, medium-, or
long-tailed according to whether $\alpha <1$, $\alpha =1$, or $\alpha >1$. When
$\alpha=1$ the relation indicated by (\ref{denQ}) may be written, in
many cases, more precisely as 
\begin{equation}\label{DQrefine}fQ(u)%
\backsim(1-u)\{\ln (1-u)^{-1}\}^{1-\beta }, 0\leq \beta \leq
1\end{equation}
 where $\beta$ is a shape
parameter. Let $L$ denote a slowly varying function from the left at 1. The precise statement associated with equation
 (\ref{denQ}) is provided by \begin{equation}\label{denQeq}
 fQ(u)=L(u)(1-u)^{\alpha}.\end{equation} Relationship (\ref{denQ}) is motivated
 by results of Andrews (1973)\nocite{Andrews}  for approximating
 the area under the tail of the distribution $F$. 
It turns out that this density-quantile representation applies for
many common distributions but not all. An example from Parzen (1979)
\nocite{Parzen} of a distribution which does not have this
density-quantile representation is $1-F(x)=\exp (-x-.75\sin x).$  To ensure that 
(\ref{denQeq}) holds, one must restrict attention to tail-monotone densities as 
discussed by Parzen (1979)\nocite{Parzen}. In addition,
Parzen (1979)\nocite{Parzen} states that the lognormal
distribution is an example with $\alpha=1$ but an expression for its
density-quantile function similar to (\ref{DQrefine}) is not
possible. Table 1 (see Parzen (1979)) gives the density-quantile function and
classification for many common distributions.
Although the classification scheme defined through (\ref{denQeq}) is of
theoretical interest, a major drawback for our purpose is that it
does not lend itself for a direct use in classifying distributions
based on data. It is difficult to check the technical assumptions
needed for (\ref{denQeq}) to hold and the various issues associated with
estimating a density arise here as well. 
In some cases, it is possible 
to estimate the tail exponent $\alpha$ in
(\ref{denQ}) under a restricted form
for $fQ(u)$. For instance, under the assumption that $fQ(u)\sim
\gamma^{-1}(1-u)^{1+\gamma}$ for $\gamma>0$, one can
estimate $\gamma$ using, for instance, the commonly used Hill estimator. As discussed earlier, however, this approach is not without problems in the case that in fact $fQ(u)\sim
\gamma^{-1}(1-u)^{1+\gamma}$ for $\gamma>0$.

\begin{table}[h]
\label{table:table3_1}
\caption[Density-quantile Functions for Various Common Distributions]{
   \renewcommand{\baselinestretch}{1} \small\normalsize
Density-quantile Functions for Various Common Distributions}

\begin{center} 

\begin{tabular}{ccc}  
Distribution        &   Density-quantile function, $fQ(u)$     &  Classification \\ 
Uniform(0,1)        & 1                                        & Short \\ 
Exponential($\gamma$)& $\frac{1}{\gamma}(1-u)$                   & Medium \\    
Logistic            & $u(1-u)$                                 & Medium \\    
Weibull($\gamma$) &$\gamma(1-u)\{log\frac{1}{(1-u)}\}^{1-\frac{1}{\gamma}}$& Medium \\   
Extreme Value       & $(1-u)log\frac{1}{(1-u)}$  & Medium \\ 
Normal              &
$\frac{1}{\sqrt{2\pi}}exp\{-\frac{1}{2}(\Phi^{-1}(u))^{2}\}\sim(1-u)(2log\frac{1}{1-
u})^{\frac{1}{2}}$  & Medium\\
Cauchy		    & $\frac{1}{\pi}\sin^{2}\pi u \sim (1-u)^{2}$ & Long\\
Pareto$(\gamma)$     & $\frac{1}{\gamma}(1-u)^{1+\gamma}$            & Long\\
Burrr$(\gamma,\tau)$&$ (1-u)^{1+\frac{1}{\gamma\tau}}$& Long
\end{tabular}
\end{center}
%\vspace{-.2in}
\end{table} 

However, as the classification based on (\ref{denQ}) yields the same
results as those
obtained using the ETT,  when the necessary technical conditions apply
for both schemes (see Parzen (1980)\nocite{Parzen80}), methodologies
based on the asymptotic distribution of the standardized maximum can
be used to classify distributions based on a random sample. In this case, however, the 
challenge is to come up with the correct sequences of constants to standardize the 
maximum. As $F$ is unknown, these will have to be estimated from the data thus
complicating the analyses. Section 4 discusses results that circumvent many
of these technical problems. The resulting methodology is easy to implement and possesses  
good operating characteristics.
Another possible drawback of a classification scheme based on (\ref{denQeq}) is that the "medium-tailed" category may be too large as discussed next.

%%%%%%%%=======================SECTION 3========================================

\section{\large Refinements and a Classification based on Extreme Spacings}

%%%%%%%%===============================================================

As the classification scheme based on (\ref{denQ}) is not sensitive enough to
distinguish, for example, among the normal, exponential, and lognormal
distributions, Schuster (1984)\nocite{Sch84}
refined Parzen's density-quantile approach for the medium tail class. This refinement
 classifies distributions such as the normal, exponential, and
lognormal into separate categories.
The following definitions using the limiting value of the failure rate function 
$r_F(x)=f(x)/\{1-F(x)\}$
give the following Refined-Parzen (RP) Classification.  Let 
\begin{equation}
\alpha =\underset{u\rightarrow 1^{-}}{\lim }-(1-u)f^{\prime
}Q(u)/[fQ(u)]^{2}, \mbox{                 and      }
\end{equation}
\begin{equation}\label{h1}
c=\lim_{u\rightarrow
1^{-}}(1-u)/fQ(u)=\lim_{u\rightarrow 1^{-}}1/r_F(Q(u)).
\end{equation}

A distribution belongs to one of the following  categories when
the given conditions hold:
\begin{center}
\begin{tabular}{ccc}
Short & $\alpha <1$ &  \\ 
Medium-Short & $\alpha =1$ & $c=0$ \\ 
Medium-Medium & $\alpha =1$ & $0<c<\infty $ \\ 
Medium-Long & $\alpha =1$ & $c=\infty $ \\ 
Long & $\alpha >1$. & 
\end{tabular}
\end{center}

The RP method classifies the normal,
exponential, and lognormal distributions as medium-short,
medium-medium, and medium-long respectively. Unfortunately, as with
Parzen's classification scheme, the RP classification cannot be
implemented easily to classify distributions from data, as it
requires estimating $fQ(u)$ and
$f'Q(u)$ for values of $u$ close to one.  

Schuster (1984)\nocite{Sch84} provided a scheme to
classify distributions by tail behavior through the asymptotic behavior of the extreme spacing (ES). That is, the difference
between the maximum and second largest data point. When the quantile
function $Q$ is differentiable in an open left interval of $1$
and if $c$ defined by (\ref{h1}) exists, Schuster categorizes
distributions by the ES as follows. \\

%%%%%%====================== Theorem 1 ==========================

\newtheorem{ESschuster}{Theorem}
\begin{ESschuster} 
\label{thm:ESschuster} 
Let $X_{1},X_{2},...,X_{n}$ be a random sample from the distribution
$F(x)$. Define  $S_{n}=X_{(n)}-X_{(n-1)}$, and assume that $c$
defined by equation
(\ref{h1}) exists. Then,\\
$
\begin{array}{cccc}
(i)\qquad c=0  & & $if and only if$ &  S_{n}=o_{p}(1),\\

(ii)\qquad c=a, & 0<a<\infty   & $if and only if$ & S_{n}=O_{p}(1),\ S_{n}\neq o_{p}(1),\\

(iii)\qquad c=\infty  & & $if and only if$  &  S_{n}\overset{p}{\rightarrow }\infty,\\
\\
\end{array}
$

\noindent where $o_{p}(1)$ denotes the sequence of random variables converges to zero in
probability, and $O_{p}(1)$ means that the sequence is bounded in probability.
\end{ESschuster}

The distribution $F$ is then said to be $\it ES\; short, ES\;medium, or\;ES\;long$, when $\it (i), (ii)$, or $\it(iii)$ hold respectively.
Theorem 1 makes the connection between the behavior of the extreme spacing and the limiting behavior of the failure rate function when $c$ defined by (\ref{h1}) exists. The failure rate goes to zero, (e.g. Pareto distribution $\bar{F}(x)=x^{-\alpha},$ with $ \alpha>0$), if and only if  the extreme spacing $S_n$ converges to infinity in probability, and the failure rate goes to  infinity (e.g. $\bar{F}(x)=e^{-e^x}$) if and only if the extreme spacing goes to zero in probability. Otherwise, the failure rate converges
to a finite positive value (e.g. $\bar{F}(x)=e^{-x})$ if and only if the extreme spacing does not converge to zero but remains bounded in probability.  

Schuster (1984)\nocite{Sch84} also made a connection between the RP
classification and the ES classification method. Using the density quantile representation
(\ref{denQeq}) and properties of slowly varying functions, it follows that 
\begin{equation}\label{aaseneta}\lim_{u\rightarrow
1^{-}}\frac{fQ(u)}{1-u}=\lim_{u\rightarrow 1^{-}}L(u)(1-u)^{\alpha
-1}= \left\{\begin{array}{lll}0 \qquad $if$ \  \alpha >1, \\
0\leq \lim_{u\rightarrow 1^{-}}L(u)\leq \infty \qquad $if$\  \alpha=1,\\ 
\infty,
\qquad $if$ \ \alpha <1.\end{array}\right.\end{equation}
Therefore the $\it ES\;short$ category consists of the RP
short and RP medium-short categories; the $\it ES\;medium$ category corresponds to
the RP medium-medium class; and the $\it ES\;long$ category consists of the RP
medium-long and RP long categories.

Thus, equation (\ref{aaseneta}) provides a simple 
 intuitive interpretation of the ES
classification. As long as the appropriate assumptions are upheld, a distribution is

\begin{center}
ES short if $1-F(x)\rightarrow 0$ faster than $f(x)\rightarrow 0$,

ES medium if $1-F(x)\rightarrow 0$ at the same rate as $f(x)\rightarrow 0$,

ES long if $1-F(x)\rightarrow 0$ slower than $f(x)\rightarrow 0$.
\end{center}

The Weibull distribution provides an example where depending on the value of the shape parameter, the Weibull distribution may be short-, medium-, or long-tailed.

\bigskip

%%%%%%=========================== EXAMPLES ============================

\noindent \textbf{Example 1:} For the Weibull distribution $\bar{F}(x)=e^{-x^{\gamma}}$,
\begin{equation}\underset{u\rightarrow 1^{-}}{\lim }\frac{1-u}{fQ(u)}=\gamma ^{-1}[-\ln
(1-u)]^{\frac{1}{\gamma }-1}=\left\{%
\begin{array}{cc}
\infty,  & 0<\gamma <1, \\ 
1, & \gamma =1, \\ 
0, & \gamma >1.%
\end{array}\right.\end{equation}
Therefore the Weibull distribution is ES short for $\gamma>1,$ ES
medium for $\gamma=1,$ and ES long for $\gamma <1.$

Additional examples of distributions classified according to the  asymptotic behavior of the extreme spacing are given in Table 2. 

%%%%%============================== END EXAMPLES =========================
There is still some degree of lack of precision in Theorem 1, and the information provided by (\ref{aaseneta}), in the sense that the case $(ii)$ in Theorem 1 and the case of $\alpha=1$ in (\ref{aaseneta}) includes medium-. short-, and long-tailed distributions. The connection between the asymptotic behavior of the failure rate and tail-heaviness of the distribution $F$ will be made precise in Lemma 5 below.

The ES classification method suggests the possibility of utilizing the asymptotic
behavior of $S_n$ to differentiate among short, medium, and long-tailed
distributions, but more specific results on the asymptotic distribution of $S_n$
are needed. 

%%%%%%===================SECTION 4=================================================

\section{\large Tail Classification using the Residual Lifetime Distribution}

%%%%%%====================================================================

Rojo (1996)\nocite{Rojo96} proposed a classification scheme based on the
asymptotic behavior of the residual life distribution. This approach circumvents many of the
technical assumptions required by previous approaches, and provides a more precise characterization 
of class membership. 

{\newtheorem{ESRojo}[ESschuster]{Definition}
\begin{ESRojo}
\label{def:ESRojo}
Define $h(t)=\underset{x\rightarrow \infty }{\lim }H_{x}(t)=\underset%
{x\rightarrow \infty }{\lim }\frac{\overline{F}(t+x)}{\overline{F}(x)}%
,\ t>0, $ when the limit exists$.$ The distribution function $F$ is considered short-tailed if $h(t)=0,$
medium-tailed if $0<h(t)<1,$ and long-tailed if $h(t)=1,$ for all $t.$\\
\end{ESRojo}
}

Although the limit $h(t)$ given in Definition \ref{def:ESRojo} exists for a fairly large class
of distribution functions, this is not always the case. The limit does not exist, for example, when there is an oscillatory behavior 
in the tail of the residual life density. Examples of distributions for
which $h(t)$ does not exist include: $\bar{F}(x)=exp(-x -.75*$sin$(x))$ and 
$\bar{F}(x)=c(1+(1+x)^{-1/2} + $sin$((1+x)^{1/2}))e^{-x}, x>0$, where $c=(2+$sin$(1))^{-1}$.

The following
results are consequences of Definition \ref{def:ESRojo}.
 Theorem \ref{thm:ESRojoCor} below combines Theorem 4.1 and
Corollary 4.2 in Rojo (1996)\nocite{Rojo96}. Hereafter, $Exp(\theta)$ will
denote the exponential distribution with parameter $\theta,$ i.e. mean $\frac{1}{\theta}.$\\

%%%%%=====================THEOREM 2=========================================

\newtheorem{ESRojoCor}[ESschuster]{THEOREM}
\begin{ESRojoCor}
\label{thm:ESRojoCor} Let $F(x)$ be a distribution
function and let $h(t)$ be as in Definition \ref{def:ESRojo}. Then,
\begin{center}
$F$ is short-tailed\qquad $ \Longleftrightarrow \qquad
S_{n}=X_{(n)}-X_{(n-1)}\overset{a.s}{\rightarrow }0,$

$F$ is medium-tailed\qquad $ \Longleftrightarrow \qquad
S_{n}=X_{(n)}-X_{(n-1)}\overset{a.s}{\rightarrow }Exp(\theta ),$

$F$ is long-tailed\qquad $  \Longleftrightarrow  \qquad
S_{n}=X_{(n)}-X_{(n-1)}\overset{a.s}{\rightarrow }\infty .$
\end{center}
\end{ESRojoCor}

Note that the results of
Theorem \ref{thm:ESRojoCor} provide a more precise characterization of 
the various classes of distributions which are in agreement with a classification
based on the asymptotic behavior of the extreme spacing $S_n$. More importantly,
 Theorem \ref{thm:ESRojoCor} delineates the asymptotic distribution of the ES for
the medium class. This is precisely the result that will lead to a methodology for 
testing the hypothesis of medium tail against either short- or long-tails in the 
next section. The asymptotic 
distribution of $S_n$ for medium-tailed distributions is perhaps not surprising since for the baseline
medium distribution, the exponential distribution, $S_n$ has an exponential distribution, for every $n$, 
with the same  parameter as the underlying distribution $F$. (see, e.g. Barlow and Proschan (1996)).
\nocite{BarlowProschan}\\

%%%%%%===========================EXAMPLES=================================
%\bigskip
\noindent \textbf{Example 1 (cont):}
For the Weibull$(\gamma)$ distribution  \begin{equation}\label{WeibClass}h(t)=\lim_{x\rightarrow \infty} 
\frac{\overline{F}(t+x)}{\overline{F}(x)}=
%\frac{e^{-(x+t)^\gamma}}{e^{-x^\gamma}}=
\lim_{x\rightarrow \infty}e^{-(x+t)^{\gamma}+x^{\gamma}}=\left\{%
\begin{array}{cc}
1  & 0<\gamma <1, \\ 
e^{-t} & \gamma =1, \\ 
0 & \gamma >1.%
\end{array}\right. \end{equation}
Therefore, for $\gamma=1$, the
Weibull(1) distribution is medium-tailed. 

For $0<\gamma<1,$ it is
 long-tailed, and for $\gamma >1$, the Weibull distribution is
short-tailed.

\bigskip

\noindent \textbf{Example 2:} The Pareto$(\gamma)$
distribution has \begin{equation}h(t)=\lim_{x\rightarrow \infty} 
\frac{\overline{F}(t+x)}{\overline{F}(x)}=\lim_{x\rightarrow \infty} \frac{(x+t)^{\gamma}}{x^{\gamma}}=
%\lim_{x\rightarrow \infty} (1+\frac{t}{x})^{\gamma}=
1\end{equation}
for all $t>0.$ Therefore the Pareto distribution is ES long by
Definition \ref{def:ESRojo}.

It is possible to refine the classification given in Definition \ref{def:ESRojo} by subdividing 
the short- and long-tailed distributions into three subclasses. This can be done
by considering, instead, the asymptotic behavior of $M(x)=\bar{F}^{-1}(e^{\bar{F}(x)})$ in the short-tailed case and the behavior of $N(x)=\bar{F}^{-1}(-1/\ln \bar{F}(x))$ in the long-tailed case. Table 2, as given in Rojo (1996), classifies several 
common distributions using the various schemes discussed so far. 

\begin{table}[h]
\label{table:table3_2}
\caption[ Refined Parzen (RP), Extreme Spacing (ES), and
Rojo's Classifications of  Distributions by Tail Behavior]{
   \renewcommand{\baselinestretch}{1} \small\normalsize
Refined Parzen (RP), Extreme Spacing (ES), and Rojo's
   Classifications of \center {Distributions by Tail Behavior}}
\begin{center}
\begin{tabular}{cccc}
\textit{Distribution} & \textit{RP} & \textit{ES} & \textit{Rojo}\\ 
Exponential & Medium-\underline{Medium} &Medium&Medium\\ 
Normal & Medium-\underline{Short} &Short&Weakly-\underline{Short}\\ 
Lognormal & Medium-\underline{Long}& Long&Weakly-\underline{Long} \\ 
Uniform & Short &Short&Super-\underline{Short} \\ 
Cauchy & Long &Long&Weakly-\underline{Long} \\ 
Extreme Value & Medium-\underline{Short} &Short&Moderately-\underline{Short} \\ 
Pareto $(\alpha <1)$ & Long & Long&Super-\underline{Long}\\  
Pareto $(\alpha =1)$ & Long & Long&Moderately-\underline{Long}\\  
Pareto $(\alpha >1)$ & Long & Long&Weakly-\underline{Long}\\  
Weibull $(\alpha <1)$&Medium-\underline{Long} &Long&Weakly-\underline{Long}\\
Weibull $(\alpha =1)$&Medium-\underline{Medium} & Medium&Medium\\
Weibull $(\alpha >1)$&Medium-\underline{Short} & Short&Weakly-Short\\
Logistic & Medium-\underline{Medium} & Medium&Medium\\
Standard Extreme Value&Medium-\underline{Short}&Short&Moderately-\underline{Short}
\end{tabular}
\end{center}
\vspace{-.2in}
\end{table}
Note that the classification scheme based on the asymptotic 
behavior of the extreme spacings, and consequently the residual life
function, is location and scale invariant. Henceforth, $\it short-, medium-$, and $\it long-tail$
will mean tail-heaviness in the sense of Theorem 3.
%%%%%%======================THE METHODOLOGY SECTION 5======================

\section{\large Testing for an ES medium tail}

%%%%%%%%%%%%%%%%%%%%%%%%%%%%%%%%%%%%%%%%%%%%%%%%%%%%%
Let $X_1, X_2, \dots, X_n$ represent a random sample from the distribution function $F$, 
and let $\bar{F}_n$ denote the empirical survival function. Consider the test statistic
\begin{equation}
\label{Test}
T_n=-\frac{\ln \overline{F}_{n}(\ln X_{(n)})(X_{(n)}-X_{(n-1)})}{\ln X_{(n)}}.
\end{equation}

This section
examines the operating characteristics of this statistic in the context of testing
hypotheses about the tail behavior of $F$.

The intuition guiding the choice of ({\ref{Test}}) as a test statistic arises 
from the following argument starting with a result of Rojo (1996)\nocite{Rojo96}. Since for a
 medium-tailed distribution $X_{(n)}-X_{(n-1)}\rightarrow$ Exp($\theta$), with probability one, for 
 some  $\theta > 0$, with $\theta$ unknown, the need arises to estimate $\theta$ to construct a test 
 statistic whose asymptotic distribution does not depend on the unknown $\theta$. Now note that,
 see Rojo (1996)\nocite{Rojo96}, for a medium-tailed distribution,
\begin{equation}\label{Rojoresult}\overline{F}(\ln
y)=y^{-{\theta }}l(y), \ \theta>0,\end{equation} where $l(y)$ is some (unknown) slowly
varying function. Therefore

\begin{equation}
\label{consistent}
\frac{-\ln \overline{F}(\ln y)}{\ln y}=\theta-\frac{\ln l(y)}{\ln y}.
\end{equation} 
The slowly varying function $l$ becomes a nuisance here, but fortunately it 
disappears in the limit  since ln$(l(y))$/ln$(y)$ converges to zero as $y\rightarrow\infty$.
Therefore, 
\begin{equation}
\label{teta}
\theta=\lim_{y\rightarrow\infty} 
-\frac{\ln \overline{F}(\ln y)}{\ln y},\end{equation} 

\hspace{-.28in} and  $\theta$ may be estimated as follows
 \begin{equation}
 \label{tetahat}
\hat{\theta}_n =-\frac{\ln \overline{F}_{n}(\ln X_{(n)})}{\ln X_{(n)}}.
\end{equation} 

The consistency of the estimator $\hat{\theta}_n$ is an immediate consequence of
the following theorem.

\newtheorem{Qconsist}[ESschuster]{Theorem}

\begin{Qconsist}
\label{thm:Qconsist}
Suppose that $\bar{F}(y)/(\bar{F}(\ln y))^{\delta} \rightarrow 0$ as $y\rightarrow \infty$ for some
$\delta>2$. Then, 
\begin{equation}
\frac{-\ln \bar{F}_n(\ln X_{(n)})}{-\ln \bar{F}(\ln X_{(n)})}\overset{P}\rightarrow 1.
\end{equation}

\end{Qconsist}

The consistency of $\hat{\theta}_n$ follows from Theorem 4 after multiplying and dividing expression  (\ref{tetahat}) by
$-\ln \bar{F}(\ln X_{(n)})$ and then using (\ref{teta}). 

Since $\theta$ is the scale parameter of the limiting distribution of $X_{(n)}-X_{(n-1)}$ 
 in the medium-tail case in Theorem 3, the
test statistic given by (\ref{Test}) 
has an asymptotic exponential
distribution with mean 1.
 That is, under the assumptions of Theorem \ref{thm:Qconsist},
 \begin{equation}
 \label{TheStat}
 -\frac{\ln \overline{F}_{n}(\ln (X_{(n)})}{\ln (X_{(n)})}(X_{(n)}-X_{(n-1)})\stackrel{D}{\rightarrow} Exp(1)\end{equation} 
 under an ES
medium-tailed distribution. To use (\ref{TheStat}) as a basis for a test for the hypothesis of medium tail, it must be verified that the conditions of Theorem \ref{thm:Qconsist} hold for medium-tailed distributions. It turns out  that $\bar{F}(y)/(\bar{F}(\ln y))^{\delta} \rightarrow 0$ as $y\rightarrow \infty$ for some
$\delta>2$ holds for all medium- and for most short- and long-tailed distributions. This follows from the following lemma which assumes that $F$ has density $f$. Let the failure rate, $f(x)/\bar{F}(x)$, associated with $F$ be denoted by $r_{F}(x)$, and let $R_{-\infty}$ denote the set of rapidly varying functions so that $f\in R_{-\infty}$ if and only if $f(\lambda x)/f(x)$ converges to zero or $\infty$, as $x\rightarrow \infty$, depending on whether $\lambda > 1$ or $\lambda < 1$, respectively. The following lemma is instrumental in proving one of the main lemmas in this section.

\newtheorem{failrate}[ESschuster]{Lemma}

\begin{failrate}
\label{lem:failrate}
Suppose that F has a density f and let $r_{F}=f/\bar{F}$ denote its failure rate. Then,\\
(i) The distribution F is short-tailed if and only if $r_{F}(t) \rightarrow \infty$ as t $\rightarrow \infty$.\\
(ii) The distribution F is medium-tailed if and only if $r_{F}(t) \rightarrow \theta$ as t $\rightarrow \infty$, for some  $0<\theta<\infty$.\\
(iii) The distribution F is long-tailed if and only if $r_{F}(t) \rightarrow 0$ as t $\rightarrow \infty.$
\end{failrate}

The following Lemma shows that the condition $\bar{F}(y)/(\bar{F}(\ln(y)))^\delta \rightarrow 0$, is achieved by all medium-tailed and most long- and short-tailed distributions. 

\newtheorem{Qmedium-tail}[ESschuster]{Lemma}

\begin{Qmedium-tail}\label{lem:Qmedium-tail}
Let F be short-, medium- or long-tailed so that $\bar{F}(\ln x) \in R_{-\infty}$, $\bar{F}(\ln x)=x^{-\theta}l_{1}(x)$, or $\bar{F}(\ln x)=l_2(x)$, respectively, for some slowly varying functions $l_i(x)$,
%=a_i(x)exp\{\int_1^{x}\frac{\varepsilon_i(t)}{t} dt\}$ for 
$i=1,2$ and some $\theta>0$. \begin{itemize}
\item[(i)] If $F$ is long-tailed, then, without loss of generality, $\bar{F}(\ln x)=\exp\{\int_1^{x}\frac{\varepsilon(t)}{t} dt\}$, with $\varepsilon(t)\rightarrow 0$ as $t\rightarrow \infty$, so that $\bar{F}(\ln x)$ is a {normalized slowly varying function} and $r_F(\ln x)=-\varepsilon(x)$.
\item[(ii)] If $F$ is short-tailed, then without loss of generality, $\bar{F}(\ln x)=\exp\{\int_1^{x}\frac{z(u)}{u}du\}$, for some $z(u) \rightarrow -\infty$  and then, $r_F(\ln x)= -z(x).$
\item[(iii)] Let $F$ be short- or long-tailed with $z(x)/z(e^x)=o(x)$ or $\varepsilon(x)/\varepsilon(e^x)=o(x)$ respectively; or suppose that $F$ is medium-tailed. Then,
\end{itemize}

\begin{equation}
\label{lemma6.0}
\bar{F}(y)/(\bar{F}(\ln y))^{\delta} \rightarrow 0  \mbox{ as } y\rightarrow \infty \mbox{ for all } \delta>0. 
\end{equation}
\end{Qmedium-tail}

Recall that a distribution $G$ is in the domain of attraction of the Frechet distribution if and only if it satisfies the von-Mises condition ( $yr_{G}(y)\rightarrow \lambda$ for some $\lambda>0$), or if $G$ is $\it tail-equivalent$ to one such distribution. Therefore, if $F$ is long-tailed, the conditions that $\varepsilon(\ln x)/\varepsilon(x)=o(x)$ is satisfied, in particular, by all those distributions satisfying the von-Mises condition. Thus, for example any $F$ with $\bar{F} \in R_{-\alpha}$ with density $f$ eventually monotone satisfies (\ref{lemma6.0}). In the case of short-tailed distributions, any distribution $F$ for which $-\ln\bar{F}(x)$ is convex in a neighborhood of infinity, will satisfy the condition that $z(\ln x)/z(x)=o(x)$ since then $r_F$ is increasing in such neighborhood and the result follows since $r_F(\ln x) = -z(x)$. An example of a short-tailed distribution that satisfies the conditions on $z$ without having increasing failure rate in a neighborhood of infinity is provided by a distribution $F$ with failure rate $r_{F}(x)=x+\ln(x)*(1+\sin(x))/x^{1/5}$.

A level $\alpha$ test of the hypothesis of an ES medium-tailed distribution can, therefore, 
be introduced using the asymptotic behavior of the test statistic $T_n$, defined by (\ref{Test}), with critical
values being the $\alpha$ and $1-\alpha$ percentiles of the
$Exp(1)$ distribution. The null hypothesis to be tested is that
$F$ is ES Medium-tailed  $vs$ either of the alternative hypotheses $H_{a_{1}}: F$ is ES Short-tailed
or $H_{a_{2}}: F$ is ES Long-tailed. The decision rules with significance levels $\alpha$ are: (1) Reject $H_{0}$ in favor of $H_{a_{1}}$ if $T_n
< -\ln(1-\alpha),$ and
(2) Reject  $H_{0}$ in favor of $H_{a_{2}}$ if  $T_n >  -\ln(\alpha)$.
Otherwise do not reject $H_{o}$.

It is clear that the asymptotic levels of the tests equal $\alpha$, and, thus, the test has point-wise robust levels.  The following two theorems prove consistency against short- and long-tailed alternatives.  In the case of short-tailed alternatives we further sub-classify them, as in Rojo (1996), according to the asymptotic behavior  of the cumulative hazard function. Thus, a short-tailed distribution $F$ is said to be {\it super-, moderately-, or  weakly-short} when
$-\ln \bar{F}(\ln x)$ is {\it rapidly-, regularly-, or  slowly-varying}. The following Theorem provides the consistency of the test for short-tailed alternatives. In the case of weakly-short distributions a mild additional condition is needed to get the result.
\newtheorem{Sconsist}[ESschuster]{Theorem}
\begin{Sconsist}
\label{thm:Sconsist}
Let $F$ be short-tailed so that $-\ln\bar{F}(\ln x)=h(x)$ with $h$ a regularly varying function with index $0\le \gamma \le \infty$.  When $\gamma=0$, suppose in addition that $\frac{r_F(\ln x)}{r_F(x)} \rightarrow 0$ as $x\rightarrow \infty$. Under the assumptions of Theorem \ref{thm:Qconsist}, the test defined by the test statistic
$T_{n}$ that rejects when
$T_{n}<-\ln (1-\alpha)$ is consistent against the class of short-tailed
alternatives. 
\end{Sconsist}

\hspace{.2in}The condition that $\frac{r_F(\ln x)}{r_F(x)} \rightarrow 0$ as $x\rightarrow \infty$ when $\bar{F}$ is a ${\it weakly-short}$ distribution is rather mild, and it is satisfied, for example, by all distributions with survival functions of the form $\bar{F}(x)=exp(-x^{(1+\alpha)})$, for $\alpha>0$; $\bar{F}(x)=exp(-x\ln_k x)$, where $\ln_k$ denotes the $kth$ iterated natural log; and $\bar{F}=\exp (-\exp(x^\alpha))$, $0<\alpha<1$.

Similar results hold for long-tailed distributions as stated in the next theorem. As in the case of short-tailed distributions, a condition is imposed on the tail of $F$ and it is seen that this condition is satisfied by a large class of distributions with either regularly or slowly varying tails. The case of long-tails with $\bar{F}$ rapidly varying, e.g. 
$\bar{F}(x)=Exp(-(x)^\alpha)$ with $0<\alpha<1$, seems to also satisfy the condition as demonstrated by many examples, but we are unable to prove the result in general.

\newtheorem{Lconsist}[ESschuster]{Theorem}
\begin{Lconsist}
\label{thm:Lconsist}
Let $F$ be long-tailed and suppose that $r_F(x)$ is eventually decreasing with
\begin{equation}
\label{condition-long}
\frac{r_F(x)}{r_F(\ln x)}=o(1).
\end{equation}
 Under the assumptions of Theorem \ref{thm:Qconsist}, the test defined by the test statistic
$T_{n}$ that rejects when
$T_{n}>-\ln (\alpha)$ is consistent against this class of long-tailed
alternatives. 
\end{Lconsist}

Other examples  that satisfy the conditions of  Theorem \ref{thm:Lconsist} include: $\bar{F}(x)=e^{-(\ln x)^\alpha}$, for $0<\alpha<1$, and regularly varying functions of the form $\bar{F}(x)=\int_1^x\frac{\varepsilon(t)}{t}dt$, where $-\varepsilon(t)\rightarrow \alpha>0$ and $-\varepsilon(t)$ is decreasing. It is possible to obtain the consistency results for regularly varying tails  by replacing the conditions  Theorem \ref{thm:Lconsist} by a different condition. This is the content of the following theorem.
\newtheorem{Lcond}[ESschuster]{Theorem}
\begin{Lcond}
\label{thm:Lcond}
Let $F$ be long-tailed with $\bar{F}$ regularly varying of exponent $\alpha \ge 0$, so that $\bar{F}(x)=c(x)exp(\int_1^x\frac{\varepsilon(t)}{t}dt)\equiv L(x)$, where $\varepsilon(t)\rightarrow -\alpha$,  and suppose that $c(x)\rightarrow c>0$, with $c(x)$ nondecreasing, and $-L'/L$ eventually non-increasing. Then, when Theorem 4 holds, the test defined in the previous theorem is consistent against these long-tailed alternatives.
\end{Lcond}
Thus the test defined by $T_{n}$ is consistent against short- and long-tailed alternatives. These are asymptotic results.    The next section provides results from a simulation study that examines the power properties for finite sample sizes.

%%%%%%%%%%%%%%%%%%%%%%%%%%%%%%%%%%%%%%%%%%%%%%%

\section{\large Simulation Results}

%%%%%%%%%%%%%%%%%%%%%%%%%%%%%%%%%%%%%%%%%%%%%%%%
The previous section discussed the asymptotic properties of Type I error control and consistency of the test against long- and short-tailed distributions.  This section investigates the  type I error, and power properties of the test for finite samples from various distributions.

Table 3 gives the rejection probabilities when 
sampling from various exponential distributions as well as
the logistic distribution and the gamma distribution with scale=1 and
shape=.7. The values given are Type I errors since all
distributions are ES medium-tailed. The probabilities for each sample
size are found from
10,000 simulations of the chosen sampling distribution.

\begin{table}[t!!]
\label{ESstatTypeI}
\caption[Power for testing for short- and long-tailed distributions when sampling  from Exponential
and Logistic Distributions]{ 
\renewcommand{\baselinestretch}{1} \small\normalsize
Type I errors -- Exponential and Logistic Distributions\\}
\begin{center}
$
\begin{array}{ccccccccccc}%cc}
n & E(100)\ S & E(100)\ L & E(1)\ S & E(1)\ L & E(.01)\ S & E(.01)\ L&Lgis\ S&Lgis\ L & G(.7)\ S & G(.7)\ L \\ % & G(1.3)\ S & G(1.3)\ L \\   
10    & .0443 & .0000 & .0541 & .0805 & .5742&.1146 & .0382 &.1427 & .1057 & .1113 \\ %& .086 & .0442 \\ 
50    & .0509 & .0001 & .0485 & .0563 & .1085&.0686 & .0447 &.0730 & .0387 & .0979  \\ %& .066 & .0316 \\ 
100  & .0515 & .0017 & .0491 & .0517 & .0675&.0642 & .0477 &.0623 & .0428 & .0928  \\%& .0637 & .0292 \\ 
250  & .0532 & .0053 & .0506 & .0541 & .0530&.0590 & .0471 &.0559 & .0382 & .0908  \\%& .0601 & .0277 \\ 
500  & .0536 & .0085 & .0542 & .0517 & .0541&.0545 & .0490 &.0594 & .0384 & .0859  \\%& .0587 & .0292 \\ 
1000& .0544 & .0100 & .0501 & .0503 & .0526&.0546 & .0462 &.0568 & .0411 & .0904  \\%& .0598 & .0298 \\
2500& .0554 & .0153 & .0514 & .0507 & .0525&.0502 & .0480 &.0572 & .0417 & .087 \\%& .0609 & .0277 \\
5000& .0523 & .0179 & .0498 & .0463 & .0477&.0499 & .0428 &.0590 & .0392 & .0853 \\%& .0578 &.0292\\
10k& .0514& .0200 & .0516 & .0458 & .0526&.0549 & .0476 & .0558    & .042 & .0906 \\%& .0607 & .028 \\
20k& .0544& .0248 & .0487 & .0508 & .0498&.0459 & .0485 & .0546    & .0377 & .083  \\%& .0676 & .0257\\
\end{array}$
\end{center}
%\vspace{-.2in}
\end{table}

Table 3 shows good performance of the test statistic 
 $- \ln \overline{F}_{n}[\ln X_{(n)}](X_{(n)}-X_{(n-1)})/\ln X_{(n)}$. 
 Most of the values are close to the desirable value of
 $\alpha=.05$ except  in the $Exp(100)$ and in the $Gamma(.7)$ cases,
 
 as well  as for a few
 instances, $Exp(1)$ and $Exp(.01)$, 
 when the sample size was extremely small, $n=10.$ 

 The case of $Gamma(.7)$ illustrates the fact that the convergence of $\ln l(y)/\ln(y)$ in (\ref{consistent}) 
 can be very slow. When this happens, the estimator $-\ln \overline{F}_{n}[\ln X_{(n)}]/{\ln X_{(n)}}$, although converging to $\delta$, it will do so rather slowly and this will be reflected on the probability of Type I error as seen in Table 3.   Despite this, the
 test statistic performs very well under various ES medium-tailed distributions.   

We now turn our attention to the power of this test statistic when
sampling from various ES short- and long-tailed distributions. Besides
tracking the power of detecting an ES short- or long-tailed
distribution, it may be just as important to notice the probability of
a serious misclassification error. A serious misclassification error
is one in which an ES short-tailed distribution sample has been classified as
long-tailed or vice-versa. Thus, both ES short and long percentages
are given. Again 10,000 simulations were used for each sample size.
Table 4 gives classification probabilities for
various shifted $Pareto(\gamma)$ distributions with survival functions
$\overline{F}(x)=\frac{1}{1+x^{\gamma}},\ x>0.$ The test statistic
shows great power against $Pareto(1)$ and $Pareto(2)$ alternatives. As expected
from the results of Heyde and Kou (2004), 
 power decreases as $\gamma$ increases. 

\begin{table}[t!!]
\label{ESstatPowerPar}
\caption[Power for testing against short- and long tails when sampling from Pareto Distributions]{
  \renewcommand{\baselinestretch}{1} \small\normalsize
Power for testing against short- and long tails -- Pareto \\}
\begin{center}
$
\begin{array}{ccccccccc}
n&P(1)\ S&P(1)\ L&P(2)\ S&P(2)\ L&P(5)\ S&P(5)\ L&P(10)\ S&P(10)\ L\\   
10  & .0126 & .5341 &.0281 &.2677&.0433 & .0411 & .0492&.0027  \\ 
50  & .0047 & .8051 &.0151 &.4547&.0346 & .1678 & .0439&.0508 \\ 
100 & .0024 & .8807 &.0119 &.5364&.0289 & .2035 & .0416&.0856 \\ 
250 & .0017 & .9420 &.0069 &.6372&.0282 & .2322 & .0396&.1169 \\ 
500 & .0007 & .9673 &.0070 &.7073&.0227 & .2660 & .0327&.1324  \\ 
1000& .0003 & .9809 &.0051 &.7726&.0219 & .2990 & .0352&.1431 \\
2500& .0002 & .9921 &.0029 &.8375&.0191 & .3562 & .0305&.1589 \\
5000& .0001 & .9968 &.0033 &.8721&.0171 & .3865 & .0300&.1733 \\
10k & .0000 & .9977 &.0019 &.9030&.0164 & .4243 & .0270&.1930 \\
20k & .0000 & .9993 &.0001 &.9284&.0140 & .4571 & .0248&.2020
\end{array}$
\end{center}
\end{table}

Table 5 gives the classification probabilities for
a $Weibull(\gamma)$ distribution with survival function
$\overline{F}(x)=e^{-x^{\gamma}},\ x>0.$ As stated previously, the
Weibull is ES short-tailed for $\gamma>1$ while ES long for
$0<\gamma<1.$ The simulations show  good power against 
the $Weibull(5)$ distribution and reasonable power for the $Weibull(1/2)$
distribution. The power decreases as $\gamma$ nears 1 as expected. 

\begin{table}[t!!]
%\vspace{-.25in}
\label{ESstatPowerWeib}
\caption[Power for testing against short- and long tails when sampling from Weibull Distributions]{
  \renewcommand{\baselinestretch}{1} \small\normalsize
Power for testing against short- and long tails -- Weibull \\}
\begin{center}
$
\begin{array}{ccccccc}
n&W(5)\ S&W(5)\ L&W(2)\ S&W(2)\ L&W(1/2)\ S&W(1/2)\ L\\   
10  & .7008 & .0000 & .1678 & .0011 & .0161&.3936\\
50  & .7589 & .0000 & .1699 & .0001 & .0127&.5050 \\ 
100 & .7862 & .0000 & .1757 & .0000 & .0087&.5450 \\ 
250 & .7999 & .0000 & .1711 & .0001 & .0101&.5855 \\ 
500 & .8071 & .0000 & .1840 & .0000 & .0085&.6171  \\ 
1000& .8177 & .0000 & .1806 & .0000 & .0069&.6500 \\
2500& .8342 & .0000 & .1851 & .0000 & .0068&.6698 \\
5000& .8428 & .0000 & .1943 & .0000 & .0067&.6844 \\
10k & .8497 & .0000 & .1910 & .0000 & .0072&.7020 \\
20k & .8519 & .0000 & .2020 & .0000 & .0049&.7158
\end{array}$
\end{center}
%\vspace{-.2in}
\end{table}

Table 6 shows the power against $U(0,1)$,
extreme value, and normal distributions. For a sample size of 100 or
larger. The test is almost perfect for detecting a $U(0,1)$ sample
as ES short. The power is unfortunately not that high for the
extreme value and normal distributions. At least in both cases for
$n>10$ the percentage ES short classifications did outnumber the ES
long ones, but the percentage of simulations that were rejected as ES
medium and classified as short-tailed from a normal distribution was
approximately 10\% for a sample size as large as 5000. 
 The lack of power in this case is addressed in the
next section. 

\begin{table}[t!!]
%\vspace{-.2in}
\label{ESstatPowerSht}
\caption[Power for testing against short- and long-tails when sampling from
 ES short distributions]{
  \renewcommand{\baselinestretch}{1} \small\normalsize
Power against short- and long-tails -- sampling from
 ES short distributions\\}
\begin{center}
$
\begin{array}{ccccccc}
n&Unif(0,1)\ S&Unif(0,1)\ L&Ext Val\ S&Ext Val\ L&Normal\ S&Normal\ L\\   
10  & .3194 & 0 & .1145 & 0     & .0465&.0533\\
50  & .7766 & 0 & .1308 & .0002 & .0629&.0148 \\ 
100 & .9459 & 0 & .1467 & 0     & .0666&.0083 \\ 
250 & .9988 & 0 & .1722 & 0     & .0816&.0041 \\ 
500 & 1 & 0     & .1816 & 0     & .0846&.0029  \\ 
1000& 1 & 0     & .1914 & 0     & .0913&.0018 \\
2500& 1 & 0     & .2077 & 0     & .0952&.0018 \\
5000& 1 & 0     & .2283 & 0     & .0984&.0006 \\
10k & 1 & 0     & .2335 & 0     & .1018&.0006 \\
20k & 1 & 0     & .2401 & 0     & .1032&.0005
\end{array}$
\end{center}
%\vspace{-.2in}
\end{table}
Finally Table 7 gives the power of the test for a
few common ES long-tailed distributions. The percentage of correct
classification for the lognormal is less than desirable, slightly
better for a $t(3)$ sample, while excellent for a Cauchy sample.  
\begin{table}[t!!]
\label{ESstatPowerLng}
\caption[Power for testing against short- and long tails when sampling from
ES Long Distributions]{
  \renewcommand{\baselinestretch}{1} \small\normalsize
Power against short- and long-tails -- sampling from
ES Long distributions\\}
\begin{center}
$
\begin{array}{ccccccc}
n&Lnorm\ S&Lnorm\ L&t(3)\ S&t(3)\ L&Cauchy\ S&Cauchy\ L  \\
10  & .0426 & .1775 & .0353 & .1833 & .0154&.4640 \\
50  & .0253 & .2828 & .0265 & .2549 & .0054&.7187 \\ 
100 & .0200 & .3329 & .0242 & .3022 & .0028&.8161 \\ 
250 & .0187 & .3910 & .0165 & .3741 & .0016&.9042 \\ 
500 & .0154 & .4437 & .0159 & .4322 & .0008&.9437  \\ 
1000& .0134 & .4953 & .0135 & .4934 & .0005&.9701 \\
2500& .0089 & .5441 & .0091 & .5740 & .0002&.9855 \\
5000& .0104 & .5826 & .0087 & .6363 & .0002&.9935 \\
10k & .0074 & .6292 & .0076 & .6879 & .0000&.9935 \\
20k & .0063 & .6574 & .0058 & .7361 & .0000&.9977 
\end{array}$
\end{center}
%\vspace{-.2in}
\end{table}

\section{\large Blocking the Data for Increased Power}

The previous section introduced a test 
to distinguish among ES short-, medium-, and long-tailed 
distribution samples by tail behavior, using the Extreme Spacing.
The test  shows
good power in distinguishing significantly different tail
behaviors. But the test showed
less capability for distinguishing a lognormal sample from an exponential
sample and  a normal sample
from an exponential sample.
This section addresses the low power values seen in the previous
section when sampling from distributions with tails which do not
differ much from the exponential. The procedure of blocking the data,
finding the test statistic for each block, and combining the block
test statistics into one test  increases the power substantially.

Notice from Table 6 that the power of detecting a
ES short-tail  when sampling from a normal distribution is
approximately 10\% for $n\le20,000$. 
Blocking the
data  into $k$ separate blocks may give rise to additional
power. Each block of size approximately $m=\frac{n}{k}$ is its own independent
subsample which will produce independent values for the test statistic $T_m$.
%\[-\frac{(X_{(m)}-X_{(m-1)})\ln\overline{F}_{m}[\ln X_{(m)}]}{\ln X_{(m)}}.\]
Under the null hypothesis of an ES medium tail, the sum of the $k$ 
block statistics can be used as the overall test statistic.

 Under the null hypothesis,
the sum of the $k$ block statistics has an asymptotic $gamma(k,1)$ distribution.  Let $TS_{j}$
be the block test statistic for block $j$ where $j=1,...,k$.
The hypotheses to be tested and corresponding decision rules are then given by 
$H_{o}: F$ is ES Medium-tailed $vs$ $H_{a_{1}}: F$ is ES Short-tailed or
$H_{a_{2}}: F$ is ES Long-tailed.
The decision rules with significance level $\alpha=5\%$ is 
Reject $H_{o}$ in favor of $H_{a_{1}}$ if 
$\Sigma_{j=1}^{k}TS_{j}<qgamma(\alpha,k,1);$
Reject $H_{o}$ in favor of $H_{a_{2}}$ if 
$\Sigma_{j=1}^{k}TS_{j}>qgamma(1-\alpha,k,1);$
otherwise do not reject $H_{o}$,
where $qgamma(p,k,1)$ is the $pth$ percentile of the gamma(k,1) distribution.

Table 8 gives the Type I errors found when sampling from the Exp(1)
and logistic distributions for the sample sizes of 500, 5000, and 20000 using
various numbers of blocks. As shown in the table the suggested number
of blocks to use is somewhere between 5 and 10, otherwise too few 
points are in each block leading to large Type I errors.  For a sample
size of 5000, it appears that up to 25 blocks can be used without causing
significant Type I errors. In what follows En stands for Exp(1) sampling with sample size $n$
and Ln stands for Logistic sampling with sample size $n$.

Tables 9-11 show that for as little as 5 or 10 blocks the
power of detecting an ES short- or long-tailed sample can increase
substantially. 

\begin{table}[t!!]
%\vspace{-.2in}
\label{BlockPowerTypeI}
\caption[Type I Errors while
Blocking the Data for Sample Sizes $n=500$, $5000$, and $20000$]{
  \renewcommand{\baselinestretch}{1} \small\normalsize
Type I Errors while
Blocking the Data for Sample Sizes $n=500$ and $5000$, and $20000$}

\begin{center}
$
\begin{array}{ccccccccc}
blocks&E500S&E500L&L500S&L500L&E5000S &E5000L&L5000S&L5000L\\   
1   & .0541 & .0545 & .0490 &.0594& .0498 & .0463 &.0428 &.0590 \\
5   & .0497 & .0538 & .0365 &.0823& .0549 & .0488 &.0381 &.0639 \\
10  & .0509 & .0660 & .0246 &.1218& .0499 & .0546 &.0332 &.0734 \\ 
25  & .0393 & .1058 & .0050 &.3469& .0463 & .0558 &.0252 &.0994 \\ 
50  & .0313 & .1344 & .0048 &.8164 & .0455 & .0674 &.0122 &.1691 \\
\end{array}$
\end{center}
\begin{center}
$
\vspace{.25in}
\begin{array}{ccccc}
blocks&Exp20k\ S&Exp20k\ L&L20k\ S&L20k\ L \\   
1   & .0487 & .0508 & .0485 &.0546 \\
5   & .0481 & .0500 & .0457 &.0597 \\
10  & .0510 & .0507 & .0382 &.0692 \\ 
25  & .0498 & .0544 & .0300 &.0861 \\ 
50  & .0485 & .0555 & .0214 &.1177 \\
100 & .0437 & .0604 & .0097 &.1885 
%\vspace{-.25in}
\end{array}$
\end{center}
\vspace{-.2in}
\end{table}  
\begin{table}[t!!]
\label{BlockPower}
\caption[Power* against Short- and Long-tailed alternatives when
Blocking the Data; $n=500$]{
  \renewcommand{\baselinestretch}{1} \small\normalsize
Power* against Short- and Long-tailed alternatives when
Blocking the Data;  $n=500$}

\begin{center}
$
\begin{array}{ccccccc}
blocks&Norm\ S&Norm\ L&Ext Val\ S&Lnorm\ L&Par(5)\ L &Weib(2) \\   
1   & .0846 & .0029 & .1816 &.4437& .2660 & .1840  \\
5   & .1653 & 0     & .7722 &.7485& .4439 & .8473  \\
10  & .1887 & 0     & .9635 &.9721& .5113 & .9941  \\ 
25  & .0950 & .0133 & .9978 &.9483& .4740 & 1  \\ 
\end{array}$
\end{center}
*{\footnotesize Not shown: Power of $0$ when testing for long tails and the distribution is Extreme Value, or $Weibull(2)$;\\
 \mbox{ }Not shown: Power $< .01$ when testing for short tails and sampling from  Lognormal, or $Pareto(1)$}
\end{table}

\begin{table}[t!!]
%\vspace{-.2in}
\label{ESBlock5000}
\caption[ES Short and Long Tail Classification Percentages while
Blocking the Data for Sample Size $n=5000$]{
  \renewcommand{\baselinestretch}{1} \small\normalsize
Power* against Short- and Long-tailed alternatives while Blocking 
the Data; $n=5000$}
\begin{center}
$
\begin{array}{cccccc}
blocks&Norm\ S&Ext Val\ S&Lnorm\ L&Par(5)\ L &Weib(2)\ S\\   
1   & .0984 & .2283 &.5826& .3865 & .1943  \\
5   & .3120 & .9406 &.9287& .7041 & .8907  \\
10  & .5172 & .9994 &.9921& .8372 & .9986  \\ 
20  & .7339 &   1    &  .9996   & .9404 & 1 \\ 
50  & .8960 &   1    &  1   &  .9924 &1  
\end{array}$
\end{center}
*{\footnotesize Not shown: 0 Long Tail Classifications for Ext Value,
		$Weibull(2)$, Normal; less than
	$1\%$ Short Tail Classifications for Lognormal,
		$Pareto(1)$}
\end{table}

\begin{table}[t!!]
\label{ESBlock20000}
\caption[ES Short and Long Tail Classification Percentages while
Blocking the Data for Sample Size $n=20000$]{
  \renewcommand{\baselinestretch}{1} \small\normalsize
Power against Short- and Lon-tailed alternatives while Blocking 
the Data; $n=20000$}

\begin{center}
$
\begin{array}{cccccc}
blocks&Norm\ S&Ext Val\ S&Lnorm\ L&Par(5)\ L &Weib(2)\ S\\   
1   & .1071 & .2442 &.6573& .4607 & .2080  \\
5   & .3804 & .9771 &.9715& .8222 & .9160  \\
10  & .6552 & 1     &.9975& .9360 & .9985  \\ 
20  & .9435 &   1    &  1   & .9902   & 1 \\ 
50  & .9966 &   1    &  1   &  .9999     &1  \\
\end{array}$
\end{center}
%\vspace{-.25in}
\end{table}
For a sample size of $n=500$ blocking the data (Table 9) increases the correct
classification to a desirable value, $>90\%$, for the extreme value,
lognormal, and $Weibull(2)$ distributions. The power increases
 for the normal and $Pareto(5)$ distributions. The reason
why the percentage does not increase over 90\% is twofold.
First, the normal and $Pareto(5)$ are very similar to an exponential
in tail behavior,
and since there are few
points in each block,  the standard error of 
$\{-\ln \bar{F}_{n}[\ln X_{(n)}]/\ln X_{(n)}\}$ increases. In other words, it is  difficult
with a sample size as small as $500$ to be able to consistently
distinguish a normal or $Pareto(5)$ tail from an exponential. The selection of the 
number of blocks is driven by a trade-off between bias and power of the tests.
Table 10 does show significant power improvement for a
sample size of 5000. More improvement is shown in Table 11 for $n=20,000$.

%%%%%%%%%%%%%%%%%%%%%%%%%%%%%%%%%%%%%%%%%%%%%%%%%%%%

\section{\large Comparison with Bryson test}

%%%%%%%%%%%%%%%%%%%%%%%%%%%%%%%%%%%%%%%%%%%%%%%%%%%%

Bryson (1974) proposed a procedure to test the hypothesis of an underlying exponential distribution
against long-tailed distributions with (increasing) linear mean residual lifetime functions. Examples of these 
long-tailed distributions are the Lomax distributions. Based on invariance considerations, the Bryson test is
defined as
\begin{equation}
\label{bryson}
T^* = \frac{\overline{X} X_{(n)}}{(n-1)\overline{X}_{GA}^2},
\end{equation}
where 
\[\overline{X}_{GA} = (\Pi_{i=1}^{n}(X_i + A_n))^{1/n}\]

\noindent with $A_n=X_{(n)}/(n-1)$. It follows from (\ref{bryson}) that the asymptotic behavior of the
 test  based on $T^*$ will be affected by the asymptotic behavior of $X_{(n)}$. One drawback of  Bryson's test
 is that its asymptotic distribution is not known and the critical values have to be simulated. Bryson (1974) provides the critical values for several small sample sizes  and three levels for the test ($\alpha=.01, .05$ and $.10$). For the 
 purpose of the present problem, this means that since we do not know that the levels of the test defined
 through $T^*$ are robust (for the class of medium-tailed distributions) then it is difficult to apply the test 
 for our purposes. Nevertheless simulation work shows that the test has good power, sometimes higher power than the test proposed here, for those distributions for which it was developed. In addition, its power
 is competitive with the power of the test defined through (\ref{teststat}) below. Its main drawback, however,
 is that it may have probability of error of Type I close to 1 when the underlying distribution is the gamma or log-gamma distributions. Thus the test may reject the null hypothesis of medium-tail in favor of a long-tail when the shape parameter of the gamma distribution is larger than 1; on the other hand, it will reject the null hypothesis 
 in favor of a short-tailed distribution, with probability of error of Type I close to 1, in the case that the shape
 parameter of the gamma is smaller than 1. For the case of the log-gamma distribution, which is long-tailed, 
 Bryson's test may reject in favor of the decision of short-tail with high probability. The reason this happens is
 that, since for the gamma distribution the centering sequence to achieve a limiting distribution for $X_{(n)}$ is given by $\log n +(\alpha-1)\log\log n - \log\Gamma(\alpha)$, then depending on whether $\alpha$ is smaller or greater than 1, the test statistic will favor a short-tail or long-tail alternative. The case of the log-gamma distribution follows in a similar manner. The following Table 12 shows the simulated quantiles for the distribution of the test statistic $T^*$ under the gamma distribution for various values of the shape parameter. In all cases, the scale parameter is set to 1. Table 13 shows the simulated quantiles for the distribution of the test statistic $T^*$ under the log-gamma distribution for various values of the scale parameter. In all cases, the shape parameter is set to 1/2.
 
The quantiles in Table 12 for shape=1, are the critical values used to implement Bryson's test. It becomes evident at once, that Bryson's test will reject, with probability close to 1, the hypothesis of medium tail in favor of short tail when the (medium-tailed) underlying distribution is gamma with shape equal to 2 and scale equal to 1 for sample sizes 100 or higher since the $97.5^{th}$ quantile for the test statistic in this case is smaller than the $2.5^{th}$ quantile of the test statistic under the standard exponential distribution. Similarly, the Bryson's test would reject the null hypothesis of medium tails with high probability if the true underlying distribution is gamma with scale $=1/2$ and shape=1. This is also obvious from Table 12 since the $97.5^{th}$ quantile for the test statistic under the null hypothesis is smaller than the $2.5^{th}$ quantile of the test statistic under the gamma with scale $=1/2$ and shape=1. Similar observations hold for the case of the log-gamma distribution.

 \begin{table}[t!!]
%\vspace{-.2in}
\label{brysontable}
\caption[Quantiles for the Bryson Test for Sample Sizes $n=50$, $100$, $500$, $5000$, and $20000$
for the gamma and log-gamma Distributions]{
  \renewcommand{\baselinestretch}{1} \small\normalsize
Quantiles for the Bryson Test for Sample Sizes $n=50$, $100$, $500$, $5000$, $10000$, and $20000$ 
for the gamma Distribution: Shape=2, 1, 1/2; Scale=1}

\begin{center}
$
\begin{array}{c|cccc||cccc||cccc}
   & .025 & .05 & .95 & .975 & .025 & .05 & .95 & .975 & .025 & .05 & .95 & .975\\   
%10 	&	0.176195205 	&	0.189258959 	&	0.37738074 	&	0.39728128 \\
\hline
50 &		0.0611 	&	0.0643 	&	0.1246 	&	0.1336 &	0.1035 	&	0.1104	&	0.2371 	&	0.2546 & 0.2042 	&	0.2212 	&	0.4588 	&	0.4852 \\
100 &	0.0386 	&	0.0407 	&	0.0752	&	0.0807 &	0.0745	&	0.0791 	&	0.1624 	&	0.1749 &  0.1751 	&	0.1889 	&	0.3768 	&	0.3968\\
%250 	&	0.043040236 	&	0.045531555 	&	0.08640843 	&	0.09304591 \\
500 &	0.0114 	&	0.0119 	&	0.0194 	&	0.0205 &	0.0271 	&	0.0283 	&	0.0506 	&	0.0543 & 0.0947 	&	0.0996 	&	0.1794 	&	0.1902 \\
%1000 &	0.016290728 	&	0.016933902 	&	0.028832644 	&	0.030773227 \\
%2500 &	0.007901103 	&	0.008161726 	&	0.01313246 	&	0.01396183 \\
5000  & 	0.0016  &		0.0016  &		0.0024 	&	0.0025 &	 0.0044 	&	0.0045 	&	0.0071 	&	0.0075  & 0.0224 	&	0.0232 	&	0.0368	&	0.0390 \\
10000 & 	0.0008 &		0.0009 &  		0.0012 	&	0.0013 &	 0.0024 	&	0.0025 	&	0.0038 	&	0.0040  & 0.0133 	&	0.0137 	&	0.0212 	&	0.0224 \\
20000 & 	0.0004 &		0.0004 &		0.0006 	&	0.0007 &   0.0013 	&	0.0013 	&	0.0020 	&	0.0021  & 0.0077 	&	0.0079 	&	0.0119 	&	0.0126\\
 
\end{array}$
\end{center}
%\vspace{-.2in}
\end{table}  
 
 \begin{table}[t!!]
%\vspace{-.2in}
\label{brysontable2}
\caption[Quantiles for the Bryson Test for Sample Sizes $n=50$, $100$, $500$, $1000$,$5000$, and $20000$ 
for the log-gamma Distribution]{
  \renewcommand{\baselinestretch}{1} \small\normalsize
Quantiles for the Bryson Test for Sample Sizes $n=50$, $100$, $500$,$1000$, $5000$, $10000$, and $20000$ 
for the log-gamma Distribution: Shape=1/2; Scale=1/6, 1}

\begin{center}
$
\begin{array}{c|cccc||cccc}
   & .025 & .05 & .95 & .975 & .025 & .05 & .95 & .975 \\   
%10 	&	0.176195205 	&	0.189258959 	&	0.37738074 	&	0.39728128 \\
\hline
50 & 	0.0239 &  0.0245 & 0.0419 &  0.0461 &	0.0682 &	0.0783 &	0.6478 &	0.7257 \\
100 & 	0.0132 &	0.0136 & 0.0240 &  0.0267 &	0.0624 &	0.0722 &	0.7082 &	0.7866 \\
500 	& 	0.0033 &  0.0034 & 0.0064 &  0.0071 &	0.0558 &	0.0671 &	0.7655 &	0.8515 \\
1000 &	0.0018 &  0.0019 & 0.0035 &  0.0039 & 	0.0564 &	0.0681 &	0.7902 &	0.8745   \\
%1000 &	0.016290728 	&	0.016933902 	&	0.028832644 	&	0.030773227 \\
%2500 &	0.007901103 	&	0.008161726 	&	0.01313246 	&	0.01396183 \\
5000  & 	0.0004 & 0.0004 & 0.0009 &  0.0010  & 	0.0558 &	0.0678 &	0.8247 &	0.9056 \\
10000 & 	0.0002 & 0.0002 & 0.0005 &  0.0005  & 	0.0559 &	0.0676 &	0.8346 &	0.9188\\
20000 & 	0.0001 & 0.0001 & 0.0002 &  0.0003  & 	0.0564 &	0.0689  &	0.8407  &	0.9215\\
 
\end{array}$
\end{center}
\vspace{-.2in}
\end{table}  
 
%%%%%%%%%%%%%%%%%%%%%%%%%%%%%%%%%%%%%%%%%%%%%%%%%%%%

\section{\large Illustrations of Data Analysis}

%%%%%%%%%%%%%%%%%%%%%%%%%%%%%%%%%%%%%%%%%%%%%%%%%%%%

This section is devoted to the analysis of
three data sets illustrating the methodologies of previous
sections. The first data set is the Secura Belgian Re data set and consists of 371 automobile claims from 1988 - 2001 from numerous European insurance companies.  Each claim was at least 1.2 million Euros.  This data, adjusted for inflation is discussed in Beirlant {\it et al.} (2004)\nocite{Teugels}.  Figures 1 and 2 show the histogram and exponential q-q plot.  Since the empirical quantiles in the right tail are greater than the corresponding exponential quantiles, the right tail appears to be longer than exponential, i.e. long-tailed.  This has been confirmed by classical techniques in Beirlant {\it et al.} (2004). \nocite{Teugels}  A Pareto-type distribution was fitted to the data and the long-tailed behavior of the data was also observed in the empirical mean residual life plots.

\begin{figure} % Figure 1n
\centering
\includegraphics[width=6.5 in, height=9in]{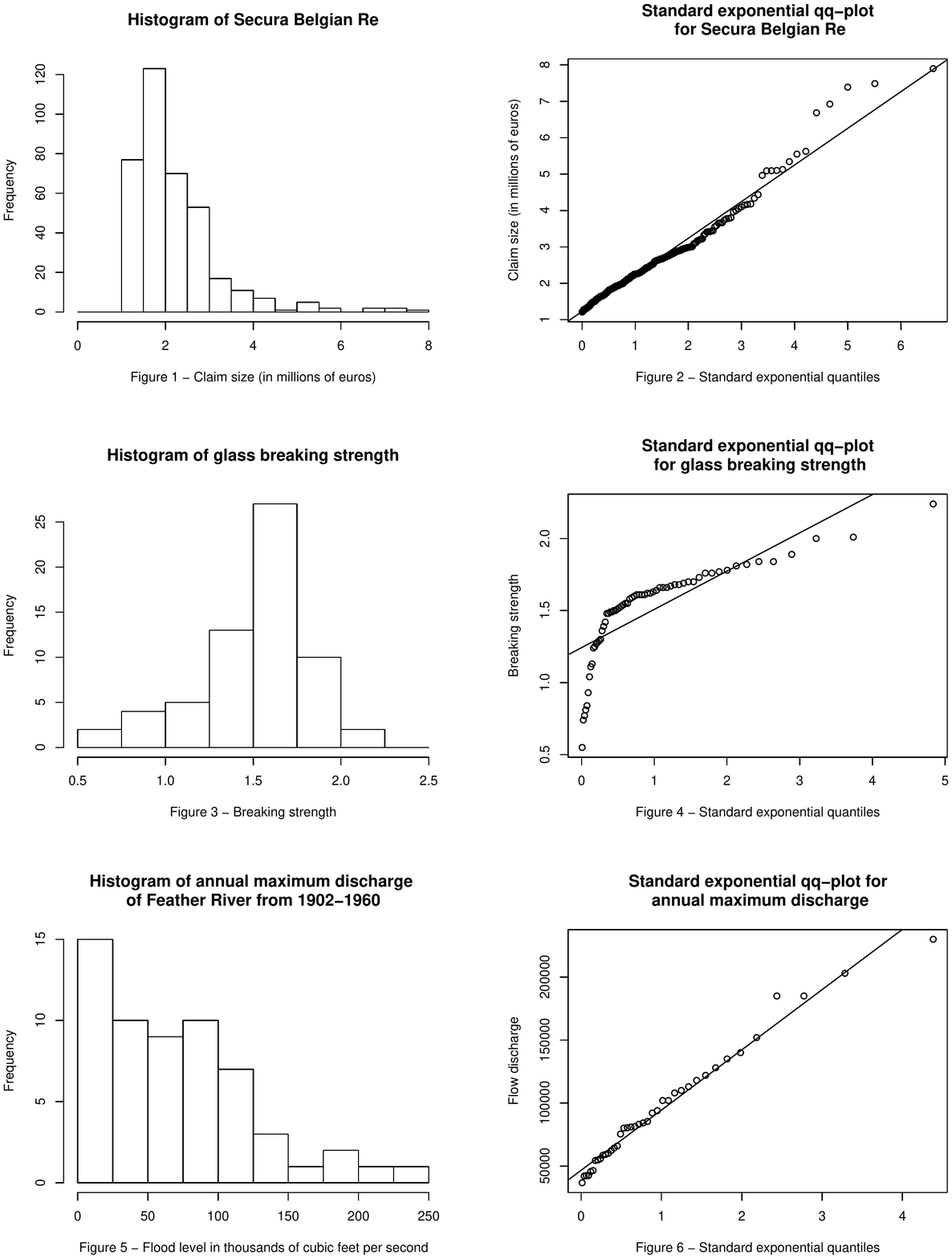} 
\label{EEM}
\end{figure}

 We consider testing for a medium tail versus a long-tail. The data, expressed in millions of Euros, was first shifted by subtracting 1.2 million. The test statistic is 
\begin{equation}
\label{teststat}
-\frac{\ln \overline{F}_{n}[\ln
X_{(n)}](X_{(n)}-X_{(n-1)})}{\ln X_{(n)}}.\end{equation} 
For claims above 1.2 million Euros the value of the test statistic is 70.40. 
Since under the null hypothesis of an ES medium-tailed the test statistic is
distributed like an $Exp(1)$, the p-value is $< .001$. Thus the distribution is classified as long-tailed.

The next example is depicted in Figure 3 and 4  that show the histogram of breaking strengths of 63 glass fibers of 1.5 cm in length.  This data appeared in Smith and Naylor (1987) \nocite{Smith} and the left tail was analyzed by Coles (2001). \nocite{Coles} Here, the the right tail is considered. The q-q plot  suggests, at first glance, a short right tail  since the empirical quantiles fall below the exponential quantiles.   

The value for test statistic is .014 which under the null hypothesis of a medium right tail gives a p-value of .014.  Therefore the test rejects the null hypothesis of a medium tail in favor of the alternative  of a short right tail.

The third example analyzes the annual maximum discharge, in thousands of cubic feet, of the Feather River from 1902 to 1960.  This data set has been described and analyzed with classical extreme value methods by Reiss and Thomas (2000). \nocite{Reiss} A Gumbel distribution was fitted to the data. The histogram and q-q plot are shown in Figures 5 and 6.  

The q-q plot suggests a medium  right tail.   After subtracting the smallest observation from the data, the test statistic yields a value of .35 with a corresponding p-value of .7 therefore not rejecting the null hypothesis of medium-tail.

\section{\large Appendix}

\vspace*{.01in}\noindent {\\

 \bf Proof of Theorem \ref{thm:Qconsist}}: Suppose that 
$\bar{F}(y)/\left(\bar{F}(\ln y)\right)^\delta \rightarrow 0$ for some $\delta>$ 2, and 
let $0 < \gamma < \frac{1}{2}$, $\varepsilon_n = n^{\gamma-\frac{1}{2}}$
and define $A_n = \{\|F_n - F\| \geq \varepsilon_n \}$ and
$B_n=\{\bar{F}(\ln X_n)\leq\varepsilon_n\}$, where $\|F_n - F\|=\sup_x|F_n(x)-F(x)|$.
Let $Z_n=\frac{-\ln\bar{F}_n(\ln X_{(n)})}{-\ln\bar{F}(\ln X_{(n)})}$ and consider first, for $\varepsilon>0$,
\begin{eqnarray*}
P\left(Z_n > 1 +\varepsilon \right) & = & P\left(Z_n > 1 + \varepsilon, A_n \right)
 +  P\left(Z_n  > 1 + \varepsilon, A_n^c\right)\\ 
 &  \leq & 2e^{-2n\varepsilon_n^2} +
P\left(Z_n> 1 + \varepsilon, A_n^c\right)
\end{eqnarray*}
Thus, it is enough to show that the second term in the last
expression goes to zero as $n \rightarrow \infty$.
Now, the second term in the last expression may be bounded from above as follows:
\begin{eqnarray*}
P\left((Z_n > 1 + \varepsilon)\cap A_n^c \cap B_n\right) & + &  P\left((\frac{-\ln(\bar{F}(\ln X_{(n)}) - \varepsilon_n)}{-\ln \bar{F}(\ln X_{(n)})}> 1 + \varepsilon)\cap A_n^c \cap B_n^c\right)\\
\vspace{3mm}
 & \leq & P(B_n) + P\left((\frac{-\ln(\bar{F}(\ln X_{(n)}) - \varepsilon_n)}{-\ln \bar{F}(\ln X_{(n)})} > 1 + \varepsilon) \cap A_n^c \cap B_n^c\right).
 \end{eqnarray*}
\begin{eqnarray*}
\vspace{-.1in}
\hspace*{-.45in} \mbox{ Now, note that } P(B_n) =   P\left(\ln X_{(n)})\geq \bar{F}^{-1}(\varepsilon_n) \right)
&  =  &  1 -  P\left(X_{(n)} < \exp\{\bar{F}^{-1}(n^{\gamma-\frac{1}{2}})\}\right)\\
&  =  &  1 - \{1-\bar{F}(\exp\{\bar{F}^{-1}(n^{\gamma-\frac{1}{2}})\})\}^{n}.
\end{eqnarray*}
But,
$n\bar{F}\left(\exp\{\bar{F}^{-1}(n^{\gamma-\frac{1}{2}})\}\right) =
o(1)$ as a consequence of the assumption
that $\bar{F}(y)/\left(\bar{F}(\ln y)\right)^\delta$ $\rightarrow 0$ as $y \rightarrow \infty$, for some $\delta>2$. To see this, set $\kappa=\gamma - 1/2$ and $u=\bar{F}^{-1}(n^\kappa)$, so that $u\rightarrow \infty$ as $n\rightarrow \infty$, and $n\bar{F}(\exp(\bar{F}^{-1}(n^\kappa)))=(\bar{F}(u))^{1/\kappa}\bar{F}(e^u)$. Finally setting $y=e^u$, then the last expression is seen to be equal to $\bar{F}(y)/(\bar{F}(\ln y))^{-1/\kappa}\rightarrow 0$ as $y\rightarrow \infty$ since $-1/\kappa>2$. 
Therefore, $P(B_n)\rightarrow 0$ as $n\rightarrow \infty$. It remains to prove that 
\[P\left((\frac{-\ln(\bar{F}(\ln X_{(n)}) - \varepsilon_n)}{-\ln \bar{F}(\ln X_{(n)})} > 1 + \varepsilon) \cap A_n^c \cap B_n^c\right) \rightarrow 0 \mbox{ as } n \rightarrow \infty.\]
\indent Write $-\ln\left(\bar{F}(\ln X_{(n)}) - \varepsilon_n\right) = -\ln
\bar{F}\left(\ln X_{(n)}\right) + \frac{\varepsilon_n}{1-\xi_n}$ where
$F\left(\ln X_{(n)}\right) < 1-\xi_n < F\left(\ln X_{(n)}\right) +\varepsilon_n$ so
that, for $0<a<1$, after setting $C_n=A_n^c\cap B_n^c$,
\begin{eqnarray*}
P\left((\frac{-\ln(\bar{F}(\ln X_{(n)}) - \varepsilon_n)}{-\ln \bar{F}(\ln X_{(n)})}
\!\! >  \!1 +  \varepsilon) \cap C_n \right)\! & = & \!P\left((1 + 
\frac{\varepsilon_n}{(1-\xi_n)(-\ln \bar{F}(\ln X_{(n)})} > 1 + \varepsilon) \cap C_n \right)\\
&\leq & P\left(\frac{\varepsilon_n}{(1-\xi_n)(-\ln\bar{F}(\ln X_{(n))})} >\varepsilon\right)\\
\end{eqnarray*}
%\begin{eqnarray*}
\[ \leq  P\left(\{\frac{\varepsilon_n}{(\bar{F}(\ln X_{(n)})-\varepsilon_n)(-\ln \bar{F}(\ln X_{(n)}))} > \varepsilon\}\cap \{n^{\frac{1}{2}-\gamma} \bar{F}(\ln X_{(n)})
> 1 + a\}\right)\]
%\end{eqnarray*}
\begin{eqnarray*}
& \hspace*{.5in}+ &  P\left(n^{\frac{1}{2} - \gamma} \bar{F}(\ln X_{(n)}) < 1 + a\right) \\
& \leq & P\left(\{\frac{1}{(n^{\frac{1}{2}-\gamma} \bar{F}(\ln X_{(n)}))(-\ln
\bar{F}(\ln X_{(n)}))}>\varepsilon\}\cap \{n^{\frac{1}{2}-\gamma} \bar{F}(\ln X_{(n)}) >
1 + a\}\right)\\
&\hspace*{.5in} + &  P\left(n^{\frac{1}{2} - \gamma} \bar{F}(\ln X_{(n)}) < 1 + a\right) \\
&\leq &  P\left(\frac{1}{a(-\ln\bar{F}(\ln X_{(n)}))} > \varepsilon\right) +
P\left(n^{\frac{1}{2} - \gamma} \bar{F}(\ln X_{(n)}) < 1 + a\right)
\end{eqnarray*}

Since $a>$ 0 while $-\ln\bar{F}(\ln X_{(n)}) \rightarrow \infty$
almost surely, the first term on the right side of the last inequality goes to zero.
For the second term, note that
\[P\left(n^{\frac{1}{2}-\gamma} \bar{F}(\ln X_{(n)}) \!<\! 1 + a\right)\! = \!
P\left(X_{(n)} \!>\! \exp\{\bar{F}^{-1}(cn^{\gamma-\frac{1}{2}})\}\right)
\! =\! 1 \!-\! \left(1 \!-\! \bar{F}(\exp\{\bar{F}^{-1}(cn^{\gamma-
\frac{1}{2}})\})\right)^n\]
where $c = 1 - a$. It then follows as before, that $\bar{F}(y)/(\bar{F}(\ln y))^\gamma \rightarrow 0$ for some $\gamma > 2$ implies that 
$n\bar{F}\left(\exp\{\bar{F}^{-1}(cn^{\gamma - \frac{1}{2}})\}\right)
\rightarrow 0$ as $n \rightarrow \infty$.  
Therefore $P\left(n^{\frac{1}{2} - \delta} \bar{F}(\ln X_{(n)}) < 1 + a
\right) \rightarrow 0$ and hence
\[P\left(\frac{-\ln \bar{F}_n (\ln X_{(n)})}{-\ln \bar{F}(\ln X_{(n)})} > 1 +
\varepsilon \right) \rightarrow 0.\]
\noindent The case  $P\left(\frac{-\ln \bar{F}_n (\ln X_{(n)})}{-\ln \bar{F}(\ln X_{(n)})} < 1 -
\varepsilon \right) =  P\left(Z_n < 1 - \varepsilon \right)$ is handled in a similar fashion.\\

Consider now
\begin{eqnarray*}
P\left(Z_n < 1 - \varepsilon \right) & = & P\left((Z_n < 1 - \varepsilon)\cap A_n \right)
+ P\left((Z_n < 1 - \varepsilon)\cap A_n^c \right)\\
& \leq & 2e^{-2n\epsilon_n^2} + P\left((\frac{-\ln (\bar{F} (\ln
X_{(n)})+\varepsilon_n)}{-\ln \bar{F}(\ln X_{(n)})} < 1 - \varepsilon)\cap A_n^c
\right)\\
& \leq & 2e^{-2n\varepsilon_n^2} + P\left(\frac{-\ln (\bar{F} (\ln
X_{(n)})+\varepsilon_n)}{-\ln \bar{F}(\ln X_{(n)})} < 1 - \varepsilon \right).
\end{eqnarray*}

As before,  write $-\ln(\bar{F}(\ln X_{(n)}) + \varepsilon_n) = -\ln(\bar{F}(\ln X_{(n)})) - \varepsilon_n/\xi_n$
where $\xi_n$ satisfies 
$\bar{F}(\ln X_{(n)})<  \xi_n < \bar{F}(\ln X_{(n)}) + \varepsilon_n.$
Then 
%\vspace{-.5in}
\begin{eqnarray*}
P\left(\frac{-\ln (\bar{F} (\ln
X_{(n)})+\varepsilon_n)}{-\ln \bar{F}(\ln X_{(n)})} < 1 - \varepsilon \right)
&=&P\left(1 - \frac{\varepsilon_n}{\xi_n(-\ln \bar{F}(\ln X_{(n)}))} < 1
- \varepsilon \right)\\
&=& P\left(\frac{\varepsilon_n}{\xi_n(-\ln
\bar{F}(\ln X_{(n)}))} > \epsilon \right)\\
&\leq& P\left(\frac{\varepsilon_n}{\bar{F}(\ln X_{(n)}) (-\ln \bar{F}(\ln
X_{(n)}))} > \varepsilon \right)\\
&=& P\left(\frac{1}{n^{\frac{1}{2}-\gamma} \bar{F}(\ln X_{(n)}) (-\ln
\bar{F}(\ln X_{(n)}))} > \varepsilon \right)\\
&\hspace*{-4.6in}\leq&  \hspace*{-2.2in}P\left(\frac{1}{-\ln \bar{F}(\ln X_{(n)})} > \varepsilon,
n^{\frac{1}{2}-\gamma}\bar{F}(\ln X_{(n)}) > 1 \right) +
P\left(n^{\frac{1}{2}-\gamma}\bar{F}(\ln X_{(n)}) < 1 \right)\\
&\hspace*{-4.6in}\leq & \hspace*{-2.2in} P\left(\frac{1}{-\ln \bar{F}(\ln X_{(n)})} > \varepsilon\right) +
P\left(n^{\frac{1}{2}-\gamma}\bar{F}(\ln X_{(n)}) < 1 \right).
\end{eqnarray*}
Similar arguments to those used before then yield the result that both terms on the right side of the above inequality  go
to zero as $n \rightarrow \infty$. Thus,
\[\frac{-\ln\bar{F}_n(\ln X_{(n)})}{-\ln\bar{F}(\ln X_{(n)})}
\overset{P}{\rightarrow} 1.\]

{\bf Proof of Lemma \ref{lem:failrate}}:
Without loss of generality it is assumed that $F$ is a life distribution. The proof follows easily by writing, after a one-step Taylor's expansion\\
\begin{equation}
\label{residlife} -\ln\{\frac{\bar{F}(t+x)}{\bar{F}(x)}\}=\int_{x}^{t+x}r_F(u)du=t*r_F(\xi),
\end{equation}
where $x<\xi<x+t$.
Consider first $(i)$. Then $\frac{\bar{F}(t+x)}{\bar{F}(x)}\rightarrow 0$, for all $t$ as $x\rightarrow \infty$.
Thus (\ref{residlife}) is equivalent to  $t*r_F(\xi)\rightarrow \infty$, as $x\rightarrow \infty$ for all $t>0$.
In the case of (ii), $\bar{F}(\ln x)=x^{-\theta}l(x)$ for some $\theta>0$ and some slowly varying function $l(x)$. Using (\ref{consistent}), it is clear that, by L'Hopital's rule,
\[\lim_{y\rightarrow \infty}r_F(y)=\lim_{y\rightarrow \infty}\frac{-\ln\bar{F}(\ln y)}{\ln y}=\theta-\lim_{y\rightarrow \infty}\frac{\ln l(y)}{\ln y}.\]
The result follows immediately since, for slowly varying $l$, $\ln l(y)/\ln y \rightarrow 0$. The converse follows immediately from (\ref{residlife}) after taking the limit as $x\rightarrow \infty$.

Case $(iii)$ also follows directly from (\ref{residlife}), since $F$ being long-tailed is equivalent to the expressions in (\ref{residlife}) converging to 0.\\

{\bf Proof of Lemma \ref{lem:Qmedium-tail}}: Consider first $(i)$. Let $\bar{F}(\ln x)=c(x)\exp{\int_1^x \frac{\varepsilon(t)}{t}dt}$, with $c(x)\rightarrow c>0$ and 
$\varepsilon(x) \rightarrow 0$, as $x\rightarrow \infty$. Since $F$ is long-tailed, $r_F(x)\rightarrow 0$ as $x\rightarrow \infty$. Therefore, 
\[\frac{r_F(\ln x)}{x}=-\frac{c'(c)}{c(x)} -\frac{\varepsilon(x)}{x}.\]
\indent This then implies that $r_F(\ln x)=\frac{-xc'(x)}{c(x)}-\varepsilon(x)\rightarrow 0$, and therefore, $xc'(x)/c(x)\rightarrow 0$. Writing $\varepsilon^*(t)=\varepsilon(t)+tc'(t)/c(t)$ it follows that
$\bar{F}(\ln x)=\exp\{\int_1^x\frac{\varepsilon^*(t)}{t}dt\}$
and the result follows since then $r_F(\ln x)=-\varepsilon^*(x)$.\\
\indent The proof of $(ii)$ is similar to that of $(i)$ after writing, for short-tailed $F$,
\[\bar{F}(\ln x)=c(x)\exp\{\int_1^x\frac{z(t)}{t}dt\}\]
with $c(x)\!\rightarrow \!c\!>\!0$ and $z(t)\!\rightarrow\! -\infty$, as $x\!\rightarrow \!\infty$, and recalling, from Lemma (\ref{lem:failrate}), that $r_F(x) \rightarrow \infty$.\\
\indent To prove $(iii)$, note that (\ref{lemma6.0}) holds if and only if 
\begin{equation}
\label{logratio}
-\ln\bar{F}(y) + \delta\ln\bar{F}(\ln y) \rightarrow \infty, \mbox{ as $y$ } \rightarrow \infty.
\end{equation}
\begin{equation}
\label{lemma6.4}
\hspace*{-.7in} \mbox{\hspace*{-.9in} Writing \qquad} -\ln\bar{F}(y) + \delta\ln\bar{F}(\ln y)=-\ln\bar{F}(y)(1-\delta(\frac{\ln\bar{F}(\ln y)}{\ln\bar{F}(y)})),\nonumber
\end{equation}
and then noticing that 
\begin{equation}
\label{failurerate}
\lim_{y\rightarrow \infty} \frac{\ln\bar{F}(\ln y)}{\ln\bar{F}(y)})=\lim_{y\rightarrow \infty} \frac{r_{F}(\ln y)}{yr_{F}(y)}=\lim_{y\rightarrow \infty}\frac{\varepsilon(y)}{y\varepsilon(e^y)}=\lim_{y\rightarrow \infty}\frac{z(y)}{yz(e^y)},
\end{equation}
where the third and fourth terms in the string of identities correspond to the cases of short- and long-tails respectively. The result for medium-tailed distributions follows immediately from Lemma \ref{lem:failrate} since in that case 
$r_F(\ln (y))/yr_F(y) \rightarrow 0$,
and for the short- and long-tailed distributions the results follow from the assumptions. \\

{\bf Proof of Theorem \ref{thm:Sconsist}}: Since $F$ is short-tailed 
$X_{(n)}-X_{(n-1)}\overset{a.s.}{\rightarrow} 0$ and 
$h(x)=-\ln \overline{F}(\ln x)$ is regularly varying with index
$\gamma$. That is,
\begin{equation}
\label{NN0} 
h(x)=x^{\gamma}l(x), \ \ 0 \le \gamma \le \infty 
\end{equation}
for some slowly varying function $l(x)$. The case $\gamma = \infty$ represents the case where $h(x)$ is rapidly varying.  It follows that 
\begin{equation}
\label{NN1}
\overline{F}^{-1}(u)=\ln h^{-1}(-\ln u)
\end{equation}
with $h^{-1}$ regularly varying with index $\frac{1}{\gamma}$.  Thus, when 
$\gamma = 0$, $h^{-1}$ is rapidly varying. Under the assumptions of Theorem \ref{thm:Qconsist},  $\frac{\ln \bar{F}_n(\ln X_{(n)})}{\ln \bar{F}(X_{(n)})} \rightarrow 1$ and it is enough to consider the behavior of

\begin{equation}
\label{NN2}
 \frac{-\ln \overline{F}(\ln X_{(n)})}{\ln X_{(n)}}(X_{(n)}-X_{(n-1)}).
\end{equation}
We prove that (\ref{NN2}) converges to zero in probablity. The cases when $\gamma=0$ or $\infty$ follow  immediately from properties of slowly and rapidly varying functions. The case of $0<\gamma<\infty$ presents the most 
technical challenges and  will be considered first. Recall that
{\it a positive function $g$ defined on some neighborhood of
$\infty$, {\bf varies smoothly with index} $\eta \in R, g \in SR_{\eta}$, if
$H(x)\overset{.}{=} \ln g(e^{x})$ $\in$ $C^{\infty}$ with
$H'(x)\rightarrow \eta$, $H^{(n)}(x)
\rightarrow 0,$ for  $n=2,3,\dots$ as $x\rightarrow \infty$.}
%\end{SVTdef}

The following theorem, (see Bingham, Goldie, and Teugels \cite {BGT}) will allow us to assume, without loss of generality, that $h(x)$ is smoothly varying, so that, as a consequence, $\lim_{t\rightarrow \infty}h'(t)t/h(t) = \gamma$.
{\bf Theorem} 
%\begin{SVT}
\label{SVT}
{\it Let $g \in R_{\eta}.$ Then there exists $g_{1},g_{2} \in SR_{\eta}$ with
$g_{1}\sim g_{2}$ \it{and} $g_{1}\leq g \leq g_{2}$ on some neighborhood of
$\infty.$  In particular, if $g \in R_{\eta},$ there exists $g^* \in
SR_{\eta}$ with $g^*\sim g.$}
%\end{SVT}
%\rm
\vspace{.1in}

Let then $0 < \gamma < \infty$ and let $U_{(1)},  \dots , U_{(n)}$ represent the order statistics from a uniform distribution on $(0,1)$. Since $X_{(n)} \overset{D}{=}\overline{F}^{-1}(1-U_{(n)})\overset{D}{=}\overline{F}^{-1}(U_{(1)})$ and
$X_{(n-1)} \overset{D}{=}\overline{F}^{-1}(1-U_{(n-1)})\overset{D}{=}\overline{F}^{-1}(U_{(2)})$, where $\overset{D}{=}$ denotes equality in distribution, expression 
(\ref{NN2}) has the same distribution as 

\begin{equation}
\label{NN3}
\frac{(\overline{F}^{-1}(U_{(1)}))^{\gamma}l(\overline{F}^{-1}(U_{(1)}))}{\ln \overline{F}^{-1}(U_{(1)})}(\overline{F}^{-1}(U_{(1)})-\overline{F}^{-1}(U_{(2)}))
\end{equation} 
where (\ref{NN3}) follows from (\ref{NN0}) and the fact that
$h(x)=-\ln \overline{F} (\ln x)$. Using a one-step Taylor's expansion, we get
\begin{eqnarray*}
\overline{F}^{-1}(U_{(1)})-\overline{F}^{-1}(U_{(2)})& = &\ln h^{-1}(-\ln U_{(1)})-\ln h^{-1}(-\ln U_{(2)})\\
 & = & \frac{U_{(2)}-U_{(1)}}{\xi_{n}h^{-1}(-\ln \xi_{n})h'(h^{-1}(-\ln \xi_{n}))}, \\
U_{(1)}<\xi_n<U_{(2)}.
\end{eqnarray*}

\hspace*{-0.25in} Therefore, (\ref{NN3}) is bounded above by 
\begin{eqnarray*}
\frac{(\overline{F}^{-1}(U_{(1)}))^{\gamma}l(\overline{F}^{-1}(U_{(1)}))}{\ln \overline{F}^{-1}(U_{(1)})}\frac{U_{(2)}-U_{(1)}}{U_{(1)}}\frac{1}{h^{-1}(-\ln \xi_{n})h'(h^{-1}(-\ln \xi_{n}))}
\end{eqnarray*}

\begin{equation}
\label{NN4}
 =\frac{(\ln h^{-1}(-\ln U_{1}))^{\gamma}l(\ln h^{-1}(-\ln U_{(1)}))}{\ln \ln h^{-1}(-\ln U_{(1)}))h^{-1}(-\ln \xi_{n})h'(h^{-1}(-\ln \xi_{n}))}\frac{U_{(2)}-U_{(1)}}{U_{(1)}}.
\end{equation}

\hspace*{-0.25in}Since 
\[\frac{U_{(2)}-U_{(1)}}{U_{(1)}}\overset{D}{=}\frac{1}{V}-1,\]
where $V\sim U(0,1)$, and since $\ln \ln h^{-1}(-\ln U_{(1)})\rightarrow\infty$ {\it a.s.}, while, as a consequence of Theorem \ref{SVT}

\begin{equation}
\label{NN5}
\frac{h'(h^{-1}(-\ln \xi_{n}))h^{-1}(-\ln \xi_{n})}{h(h^{-1}(-\ln \xi_{n}))}\rightarrow\gamma > 0
\end{equation}
then, to show that 
\[\frac{-\ln \overline{F}(\ln X_{(n)})}{\ln X_{(n)}}(X_{(n)}-X_{(n-1)})\overset{P}{\rightarrow}0\]
 it is enough to show that 
\begin{equation}
\label{NN6}
\frac{(\ln h^{-1}(-\ln U_{(1)}))^{\gamma}l(\ln h^{-1}(-\ln U_{(1)}))}{-\ln \xi_{n}}\overset{P}{\rightarrow}0.
\end{equation}

To verify (\ref{NN5}), 
write $h(x)=-\ln \overline{F}(\ln x)$ so that 
$h'(x)=\frac{r_F(\ln x)}{x}$ and
$h^{-1}(t)=\exp\{\overline{F}^{-1}(e^{-t}),\}$
it follows that, after setting $t=\bar{F}^{-1}(\xi_n)$,
\[\frac{h'(h^{-1}(-\ln \xi_{n}))h^{-1}(-\ln \xi_{n})}{-\ln \xi_{n}}=\frac{r_F(t)}{-\ln \overline{F}(t)}
=\frac{d}{dt}\ln (-\ln\overline{F}(t))|_{t=\overline{F}^{-1}(\xi_{n})}=\frac{d}{dt}\ln h(e^{t}).\]

\hspace*{-.25in}Thus,  (\ref{NN5}) follows from Theorem \ref{SVT}
with $\eta=\gamma$.
To prove (\ref{NN6}) rewrite as
\[\frac{(\ln h^{-1}(-\ln U_{(1)}))^{\gamma}l(\ln h^{-1}(-\ln U_{(1)}))}{-\ln U_{(1)}}\frac{-\ln U_{(1)}}{-\ln(\xi_{n})}\]
which is bounded above by 
\[\frac{(\ln h^{-1}(-\ln U_{(1)}))^{\gamma}l(\ln h^{-1}(-\ln U_{(1)}))}{-\ln U_{(1)}}\frac{-\ln
U_{(1)}}{-\ln U_{(2)}}.\]
Now observe that 
%\begin{equation}
%\label{NN7}
$P(\frac{-\ln U_{(1)}}{-\ln U_{(2)}}>2)=o(1)$
%\end{equation}
and in fact,  $P(\frac{-\ln U_{(1)}}{-\ln U_{(2)}}>3 \ i.o.)=0$. Therefore, writing $t_{n}=-\ln U_{(1)}$,

\[\frac{(\ln h^{-1}(t_{n}))^{\gamma}l(\ln h^{-1}(t_{n}))}{t_{n}}\overset{a.s.}{\rightarrow} 0\]
 since  $h^{-1}$ is $R_{\frac{1}{\gamma}}$ so that $\ln h^{-1}$ is
slowly varying, and hence $l(\ln h^{-1}(t_{n}))$ and
$(\ln h^{-1}(t_{n}))^{\gamma}$ are slowly varying.
It follows that (\ref{NN6}) is true since
$l(x)/x\rightarrow 0$ for slowly varying $l$.

Consider now the case of $\gamma=0$. Since %
$\frac{-\ln U_(1)}{-\ln \xi_n}\rightarrow 1$ in probability,
it follows from (\ref{NN4}) that to show that  (\ref{NN2}) converges to zero in probability, it is enough to show that
\begin{equation}
\label{next-to-last}
\frac{h(\ln h^{-1}(-\ln U_{(1)}))}{\ln\ln h^{-1}(-\ln U_{(1)})h^{-1}(-\ln U_{(1)})h'(h^{-1}(-\ln U_{(1)}))}\rightarrow 0,
\end{equation}
and this follows directly, after writing $y=\ln h^{-1}(-\ln U_{(1)})$, from the assumptions in the case of $\gamma=0$, since in this case, $r_F(x)=e^xh'(e^x)$ and $-\ln\bar{F}(\ln x)/\ln x\sim r_F(\ln x)$.

Finally, consider the case of $\gamma=\infty$. That is, suppose that $-\ln \bar{F}(\ln x)$ is rapidly varying. The condition given by (\ref{next-to-last})
is seen to be equivalent to 
\begin{equation}
\label{last}
\frac{h(x)}{e^x h'(e^x)\ln x}\rightarrow 0, \mbox{ as } x\rightarrow \infty.
\end{equation}
Recall that a rapidly varying function $h(x)$ may be written as $c(x)exp({\int_1^x \frac{z(t)}{t}dt})$, with $c(x)\rightarrow c$ and  $z(t)>0$, $z(t)\rightarrow\infty$, as $t\rightarrow\infty$. Assuming without loss of generality that $c(x)=c$, it is clear that $e^xh'(e^x)=c*z(e^x)h(e^x),$ and hence, 
\begin{equation}
\label{last-last}
\frac{h(x)}{e^x h'(e^x)\ln x}=\frac{h(x)}{z(e^x)h(e^x)\ln x}.
\end{equation}
Since $z(x)\rightarrow\infty$, as $x\rightarrow \infty$, while $h$ is non-decreasing, it follows that 
(\ref{last}) holds. \\%\vspace{-.15in} \begin{flushright} \qed \end{flushright}

{\bf Proof of Theorem \ref{thm:Lconsist}}: For long-tailed $F$, $\bar{F}(\ln x )=L(x)$ for some slowly varying $L$. Therefore, $-\ln \bar{F}(\ln x)=-\ln L(x)$ and hence $r_F(\ln x)=-xL'(x)/L(x)$,  where $L(x)=c(x)\exp(\int_1^x \frac{\varepsilon(t)}{t}dt)$. Under the conditions of Theorem \ref{thm:Qconsist}  consider $\frac{-\ln \bar{F}(\ln X_{(n)})}{\ln X_{(n)}}$ instead of $\frac{-\ln \bar{F_n}(\ln X_{(n)})}{\ln X_{(n)}}$ and write
\begin{equation}
\frac{-\ln \bar{F}(\ln X_{(n)})}{\ln X_{(n)}}=\frac{-\ln \bar{F}(\ln X_{(n-1)})}{\ln X_{(n-1)}}+(X_{(n)}-X_{(n-1)})\{\frac{r_F(\ln \xi_n)}{\xi_n\ln \xi_n}+\frac{\ln \bar{F}(\ln \xi_n)}{\xi_n(\ln \xi_n)^2}\}\nonumber
\end{equation}
for some $\xi_n$ with $X_{(n-1)}<\xi_n<X_{(n)}$. Note that, almost surely, as $n\rightarrow \infty$
\begin{equation}
\frac{\ln \bar{F}(\ln \xi_n)}{\xi_n(\ln\xi_n)^2}\sim \frac{-r_F(\ln \xi_n)}{\xi_n(2\ln \xi_n+(\ln \xi_n)^2)}.\nonumber
\end{equation}
Therefore, almost surely, for sufficiently large $n$,
%\label{nby(n-1)}
\[
T_n\ge \frac{-\ln \bar{F}(\ln X_{(n-1)})}{\ln X_{(n-1)}}(X_{(n)}-X_{(n-1)})\nonumber
\]

Since $\bar{F}^{-1}(u)=\ln(L^{-1}(u))$, we can write, using a one-step Taylor's expansion,
\begin{eqnarray*}
%\label{long}
\frac{-\ln \bar{F}(\ln X_{(n-1)})}{\ln X_{(n-1)}}(X_{(n)}-X_{(n-1)}) & = &\frac{-\ln L(X_{(n-1)})}{\ln X_{(n-1)}}(\ln L^{-1}(U_{(1)})-\ln L^{-1}(U_{(2)})) \\
& = & \frac{-\ln L(\ln L^{-1}(U_{(2)}))}{\ln \ln L^{-1}(U_{(2)})}\frac{(U_{(1)}-U_{(2)})}{L'(L^{-1}(\psi_n))L^{-1}(\psi_n)},
\end{eqnarray*}
where $U_{(1)}, U_{(2)}$ are the first and second order statistics from a $U\sim(0,1)$ random sample, and $U_{(1)}<\psi_n<U_{(2)}$. Writing $U_{(1)}-U_{(2)}=-(1-V)U_{(2)}$, where $V\sim U(0,1)$ with $V$  independent of $U_{(2)}$, and since $\psi_n<U_{(2)}$, and $L'(x)<0$, it follows that $T_n$ is at least as large as
\begin{eqnarray}
\label{long1}
\frac{-\ln L(\ln L^{-1}(U_{(2)}))}{\ln \ln L^{-1}(U_{(2)})}\frac{(1-V)LL^{-1}(\psi_n)}{-L'(L^{-1}(\psi_n))L^{-1}(\psi_n)}.
\end{eqnarray}
Since $r_F$ is eventually decreasing, $-xL'(x)/L(x)$ is eventually decreasing and since $L^{-1}(\psi_n)\ge L^{-1}(U_{(2))}$, it follows that 
\[\frac{-L^{-1}(U_{(2)})L'L^{-1}(U_{(2)})}{LL^{-1}(U_{(2)})}\ge \frac{-L^{-1}(\psi_n)L'L^{-1}(\psi_n)}{LL^{-1}(\psi_n)},\]
and hence, (\ref{long1}) implies that, almost surely, for large $n$, after setting $Y_n=L^{-1}(U_{(2)})$, 
\begin{equation}
\label{long2}
T_n\ge \frac{-\ln L(\ln Y_n)}{\ln \ln Y_n}\frac{(1-V)L(Y_n)}{-L'(Y_n)Y_n}
\end{equation}
with $V$ independent of $Y_n$, and $Y_n\rightarrow \infty$ almost surely. Now 
\[\frac{L(y)}{-L'(y)y}=\frac{1}{r_F(\ln y)}, \mbox{ while } \frac{-\ln L(\ln y)}{\ln \ln y}\sim \frac{-L'(\ln y)\ln y}{L(\ln y)}=r_F(\ln \ln y).\]
Therefore, it follows that the right side of (\ref{long2}) is asymptotically equivalent, almost surely, to
\[(1-V)\frac{r_F(\ln \ln Y_n)}{r_F(\ln Y_n)}.\]
The result then follows from the assumption that $r_F(y)=o(r_F(\ln y)$. 

{\bf Proof of Theorem \ref{thm:Lcond}}: As in the proof of the previous theorem, we can write,  almost surely, for sufficiently large $n$
\begin{equation}
T_n\ge \frac{-\ln \bar{F}(\ln X_{(n-1)})}{\ln X_{(n-1)}}(X_{(n)}-X_{(n-1)}).
\end{equation}

Since $\bar{F}^{-1}(u)=L^{-1}(u)$, we can write, using a one-step Taylor's expansion,
\begin{equation}
\label{long3}
\frac{-\ln \bar{F}(\ln X_{(n-1)})}{\ln X_{(n-1)}}(X_{(n)}-X_{(n-1)}) % & = & %\frac{-\ln L(\ln L^{-1}(U_{(2)}))}{\ln X_{(n-1)}}(\ln L^{-1}(U_{(1)})-\ln L^{-1}(U_{(2)})\nonumber \\
\end{equation}
\begin{equation}
= \frac{-\ln L(\ln L^{-1}(U_{(2)}))}{ \ln L^{-1}(U_{(2)})}\frac{(U_{(1)}-U_{(2)})}{L'(L^{-1}(\psi_n))} ,\\
\end{equation}
\begin{equation}
\label{Taylor}
 =   \frac{-\ln L(\ln L^{-1}(U_{(2)}))}{ \ln L^{-1}(U_{(2)})}\frac{(1-V)U_{(2)}}{-L'(L^{-1}(\psi_n))},
\end{equation}
where $U_{(1)}, U_{(2)}$ are the first and second order statistics from a $U\sim(0,1)$ random sample, and $U_{(1)}<\psi_n<U_{(2)},$ with $V$ independent of $U_{(2)}$.

But, since $L^{-1}$ is decreasing and, therefore, $L^{-1}(\psi_n)>L^{-1}(U_{(2)})$, and $-L'/L$ is eventually decreasing, then 
\[\frac{U_{(2)}}{-L'(L^{-1}(\psi_n)}\ge\frac{L(L^{-1}(\psi_n))}{-L'(L^{-1}(\psi_n)}\ge \frac{L(L^{-1}(U_{(2)}))}{-L'(L^{-1}U_{(2)})}. \]
It follows from (\ref{long3}) that, after setting $Y_n=L^{-1}(U_{(2)})$
\begin{eqnarray*}
\frac{-\ln \bar{F}(\ln X_{(n-1)})}{\ln X_{(n-1)}}(X_{(n)}-X_{(n-1)})& \ge & \frac{-\ln L(\ln Y_n)L(Y_n)}{-\ln Y_n L'(Y_n)}
=  \frac{-\ln L(\ln Y_n)}{-\ln Y_n (\frac {-c'(Y_n)}{c(Y_n)} -\frac{\varepsilon(Y_n)}{Y_n})}\\
& \ge &  \frac{-\ln L(\ln Y_n) Y_n}{\ln Y_n (-\varepsilon(Y_n))} \rightarrow \infty, 
\end{eqnarray*}
since $-\varepsilon(Y_n)\rightarrow \alpha \ge 0$ and $\frac{-\ln L(\ln Y_n) Y_n}{\ln Y_n} \rightarrow \infty$, almost surely. 
\\

\begin{center}
\textsc{\bf Acknowledgements}
\end{center}
%\vspace{-.2in}
Thanks to David J. Kahle for helpful discussions.\\

\small\normalsize
\vspace{-.2in}
\addcontentsline{toc}{chapter}{\numberline{}Bibliography}
%\nocite{*}  % UNCOMMENT TO INCLUDE ALL REFERENCES IN FILE 
             % AND NOT JUST THOSE CITED
\bibliographystyle{plain}
\bibliography{XPapersources}

%%%%%%%%%%%%%%%%%%
\end{document}